\documentclass[12pt]{amsart}
\usepackage{cases}
\usepackage{amsthm}
\usepackage{amsmath}
\usepackage{amscd}
\usepackage{graphicx}
\usepackage[mathscr]{eucal}
\usepackage[colorlinks,linkcolor=blue,citecolor=blue, pdfstartview=FitH]{hyperref}

\input xy
\xyoption{all} \numberwithin{equation}{section}
\numberwithin{figure}{section}
\usepackage{float}

\begin{document}

\title[BPS invariants from framed links]
{BPS invariants from framed links}

\author[Kai Wang and Shengmao Zhu]{
Kai Wang and Shengmao Zhu}

\address{Department of Mathematics \\
Zhejiang Normal University  \\
Jinhua Zhejiang,  321004, China }

\email{meetmath@zjnu.edu.cn, szhu@zju.edu.cn}

\begin{abstract} 
In this article, we investigate the BPS invariants associated with framed links. 
We extend the relationship between the algebraic curve (i.e. dual $A$-polynomial) and the BPS invariants of a knot investigated in \cite{GKS} to the case of a framed knot. With the help of the framing change formula for the dual $A$-polynomial of a framed knot, we give several explicit formulas for the extremal $A$-polynomials and the BPS invariants of framed knots.  As to the framed links, we present several numerical calculations for the Ooguri-Vafa invariants and BPS invariants for framed Whitehead links and Borromean rings and verify the integrality property for them.   
\end{abstract}

\maketitle

\theoremstyle{plain} \newtheorem{thm}{Theorem}[section] \newtheorem{theorem}[%
thm]{Theorem} \newtheorem{lemma}[thm]{Lemma} \newtheorem{corollary}[thm]{%
Corollary} \newtheorem{proposition}[thm]{Proposition} \newtheorem{conjecture}%
[thm]{Conjecture} \theoremstyle{definition}
\newtheorem{remark}[thm]{Remark}
\newtheorem{remarks}[thm]{Remarks} \newtheorem{definition}[thm]{Definition}
\newtheorem{example}[thm]{Example}





\tableofcontents

\newpage

\section{Introduction}
There are many amazing relations between quantum field theory, string theory and knot theory. Two aspects of the polynomial knot invariants (for instance, HOMFLYPT invariants, superpolynomials), the special geometry of algebraic curves associated to a knot $\mathcal{K}$ and the integrality of the BPS invariant associated to $\mathcal{K}$ are discussed in \cite{GKS}. The algebraic curve associated to a knot is a curve obtained as a classical limit of recursions of knot polynomials. Such recursions are proved to exist for colored Jones polynomials \cite{GL05,Hik04}, colored HOMFLYPT invariants \cite{GLL18}.  Moreover, for the colored HOMFLYPT invariant, the above algebraic curve is conjectured to agree with the augmentation polynomial of a knot \cite{AV12,AENV13}. On the other hand, by the BPS invariants of a knot, we mean the BPS degeneracies (or Ooguri-Vafa invariant, Labastida-Marino-Ooguri-Vafa (LMOV) invariants)  of the a knot. 

In \cite{GKS}, Garoufalidis-Kucharski-Sulkowski gave a method to derive the BPS invariants of a knot from its algebraic curve, which is determined by a polynomial, sometimes referred to as the dual $A$-polynomial. But in general it is difficult to write down an explicit formula for the BPS invariant,  so \cite{GKS} focused on the extremal BPS invariants which can be expressed in terms of coefficients of maximal or minimal powers of $a$ in colored HOMFLYPT invariant. The algebraic curve that encodes extremal BPS invariants is called the extremal $A$-polynomial. Explicit formulas for the extremal BPS invariants and the extremal $A$-polynomials of the twist knots $\mathcal{K}_p$ were presented in \cite{GKS}.

The goal of this article is to generalize the above relationship between the algebraic curves and BPS invariants to the case of a framed knot. Moreover, we also investigate the integrality conjecture for BPS invariants of framed links. More precisely, we will use the framed colored HOMFLYPT invariant which depends on the framing of a framed knot to rewrite the above story. A framed knot $\mathcal{K}^\tau$ is a knot $\mathcal{K}$ with a specific framing $\tau\in \mathbb{Z}$. Given a partition $\lambda$, let $H_\lambda(\mathcal{K};q,a)$ denote the (framing-independent) colored HOMFLYPT invariants of a knot $\mathcal{K}$. Hence $H_\lambda(\mathcal{K};q,a)=H_\lambda(\mathcal{K}^\tau;q,a)$ for any specific framing $\tau$.  In this article, we introduce the framed colored HOMFLYPT invariants of
$\mathcal{K}^\tau$ as follows
\begin{align} \label{framedknotformula}
\mathcal{H}_\lambda(\mathcal{K}^\tau,q,a)=(-1)^{|\lambda|\tau}q^{\frac{\kappa_\lambda\tau}{2}}H_{\lambda}(\mathcal{K},q,a),
\end{align}
where
$\kappa_\lambda=\sum_{i=1}^{l(\lambda)}\lambda_i(\lambda_i-2i+1)$.  
In particular, when $\lambda=(n)$, we have
\begin{align}
    \mathcal{H}_n(\mathcal{K}^\tau;q,a)=(-1)^{n\tau}q^{\frac{n(n-1)}{2}\tau}H_{n}(\mathcal{K},q,a).
\end{align}
We should stress that the factor $(-1)^{n\tau}$ is very important in the above definition of the framed colored HOMFLYPT invariants, it ensures the integrality structure of the BPS invariants.

The generating function of $\mathcal{H}_n(\mathcal{K}_\tau;q,a)$ determines the dual $A$-polynomial $\mathcal{A}_{\mathcal{K}^\tau}(x,y,a)$. In particular, by its definition, when $\tau=0$, $\mathcal{A}_{\mathcal{K}^0}(x,y,a)$ is just the dual $A$-polynomial $\mathcal{A}_\mathcal{K}(x,y,a)$ introduced in \cite{GKS}. Moreover, we  can derive the following explicit framing change formula 
\begin{align} \label{formula-Aframed}
    \mathcal{A}_{\mathcal{K}^\tau}(x,y,a)=\mathcal{A}_{\mathcal{K}^0}((-1)^\tau y^{2\tau} x,y,a). 
\end{align}
Therefore,  the BPS invariants for the framed knot $\mathcal{K}^\tau$ can be read from the formula (\ref{formula-Aframed}), as well as the extremal BPS invariants.

We show that such a generalization is nontrivial. First, we compute the dual $A$-polynomial and BPS invariants for the framed unknot $U^\tau$ in detail. Although the BPS invariant for unknot $U$ is trivial \cite{GKS}, it is highly nontrivial for the framed unknot $U^\tau$. 

\begin{proposition}
(1)   The dual $A$-polynomial for framed unknot $U^\tau$ is given by 
\begin{align}
  \mathcal{A}_{U^\tau}(x,y,a)=(-1)^\tau(y^2-1)-x(a^{\frac{1}{2}}y^2-a^{-\frac{1}{2}})y^{2\tau}.   
\end{align}   
(2) The BPS invariant for $U^\tau$ is given by   
\begin{align} \label{formula-brm}
b_{r,m}(U^\tau)=\frac{1}{r^2}\sum_{d|r,m}\mu(d)c_{\frac{r}{d},\frac{m}{d}}(U^\tau)    
\end{align}
where $\mu(d)$ is the M\"obius function and

\makeatletter
\let\@@@alph\@alph
\def\@alph#1{\ifcase#1\or \or $'$\or $''$\fi}\makeatother
\begin{subnumcases}
{c_{r,m}(U^\tau)=}
(-1)^{r\tau+r+\frac{r+m}{2}}\binom{r}{\frac{r+m}{2}}\binom{r\tau+\frac{r+m}{2}-1}{r-1}, &$r+m=\text{even}$, $|m|\leq r$ \label{eq:a1} 
\nonumber \\\nonumber
0, & otherwise.\label{eq:a2}
\end{subnumcases}
\makeatletter\let\@alph\@@@alph\makeatother
\end{proposition}

\begin{proposition}
(1) The extremal $A$-polynomial for framed unknot $U^\tau$ is given by 
\begin{align}
 \mathcal{A}_{U^\tau}^+(x,y)=(-1)^\tau(y^2-1)-xy^{2\tau+2},\ \  \mathcal{A}^{-}_{U^\tau}(x,y)=(-1)^{\tau}(y^2-1)+xy^{2\tau}.  
\end{align}

(2) The extremal invariant for $U^\tau$ is given by 
\begin{align} \label{formula-brpm}
b^+_{r}(U^\tau)&=\frac{1}{r^2}\sum_{d|r}\mu\left(\frac{r}{d}\right)(-1)^{d\tau}\binom{d(\tau+1)-1}{d-1},\\ \nonumber   
b^-_{r}(U^\tau)&=\frac{1}{r^2}\sum_{d|r}\mu\left(\frac{r}{d}\right)(-1)^{d(\tau+1)}\binom{d\tau-1}{d-1}.    
\end{align}
\end{proposition}
Comparing the formulas (\ref{formula-brm}) and (\ref{formula-brpm}), we have $b^\pm_{r}(U^\tau)=b_{r,\pm r}(U^\tau)$.

Then, we compute the extremal $A$-polynomials and the extremal BPS invariants for framed twist knots, extending the results in \cite{GKS}.  

\begin{proposition}
For the framed twist knot $\mathcal{K}_p^\tau$ with framing $\tau\in \mathbb{Z}$.

(1) The extremal A-polynomial of $\mathcal{K}_p^\tau$ is given by
\begin{align}
\mathcal{A}_{\mathcal{K}_p^\tau}^-(x,y)&=(-1)^\tau xy^{2\tau}-y^4+y^6, \\\nonumber \mathcal{A}_{\mathcal{K}_p^\tau}^+(x,y)&=1-y^2+(-1)^\tau xy^{4|p|+2+2\tau}
\end{align}
for $p\leq 1$ and
\begin{align}
\mathcal{A}_{\mathcal{K}_p^\tau}^-(x,y)&=1-y^2-(-1)^\tau x y^{4+2\tau}, \\\nonumber \mathcal{A}_{\mathcal{K}_p^\tau}^+(x,y)&=1-y^2-(-1)^\tau xy^{4p+4+2\tau}
\end{align}
for $p\geq 2$. 

(2) The extremal BPS invariants of $\mathcal{K}_p^\tau$  is given by 
\begin{align}
 b_{r}^-(\mathcal{K}^\tau_p)&=-\frac{1}{r^2}\sum_{d|r}\mu\left(\frac{r}{d}\right)(-1)^{d\tau}\binom{d(3-\tau)-1}{d-1}, \\\nonumber 
b_{r}^+(\mathcal{K}^\tau_p)&=\frac{1}{r^2}\sum_{d|r}\mu\left(\frac{r}{d}\right)(-1)^{d\tau}\binom{d(2|p|+1+\tau)-1}{d-1}    
\end{align}
for $p\leq -1$ and 

\begin{align}
 b_r^-(\mathcal{K}^\tau_p)&=\frac{1}{r^2}\sum_{d|r}\mu\left(\frac{r}{d}\right)(-1)^{d(\tau+1)}\binom{d(\tau+2)-1}{d-1}, \\\nonumber
 b_{r}^+(\mathcal{K}^\tau_p)&=\frac{1}{r^2}\sum_{d|r}\mu\left(\frac{r}{d}\right)(-1)^{d(\tau+1)}\binom{d(\tau+2+2p)-1}{d-1}
\end{align}
for $p\geq 2$.
\end{proposition}

Finally, based on the previous works of the second author \cite{LZhu19,Zhu19}, we have the following.
\begin{theorem}
All the BPS invariants $b_{r,m}(U^\tau)$, $b_r^{\pm}(U^\tau)$, $b^{\pm}_r(\mathcal{K}_p^\tau)$ listed above are integral.     
\end{theorem}

As another generalization, we consider the integrality of the BPS invariants for framed links. 

Given a framed link $\mathcal{L}$ with the $L$ components $\mathcal{K}_1,..,\mathcal{K}_L$. Suppose that the framing of $\mathcal{K}_\alpha$ is given by $\tau_\alpha$, we let $\vec{\tau}=(\tau_1,...,\tau_L)$ and denote the framed link $\mathcal{L}$ by $\mathcal{L}^{\vec{\tau}}$. 
The framed colored HOMFLYPT invariant (with symmetric representations) for $\mathcal{L}^{\vec{\tau}}$ is defined by 
\begin{align}
\mathcal{H}_{r_1,...,r_L}(\mathcal{L}^{\vec{\tau}};q,a)=(-1)^{\sum_{\alpha=1}^Lr_\alpha|\tau_\alpha|}q^{\sum_{\alpha=1}^L\frac{r_\alpha(r_\alpha-1)\tau_\alpha}{2}}H_{r_1,...,r_L}(\mathcal{L};q,a),  
\end{align}
where $H_{r_1,...,r_L}(\mathcal{L};q,a)$ is the framing-independent colored HOMFLYPT invariant with each component $\mathcal{K}_\alpha$ colored by symmetric representation labeled by a positive integer $r_\alpha\in \mathbb{N}$.  

We consider the following generating function for  $\mathcal{H}_{r_1,...,r_L}(\mathcal{L}^{\vec{\tau}};q,a)$, 
\begin{align} \label{formula-Zlink}
\mathcal{Z}(\mathcal{L}^{\vec{\tau}};q,a;x_{1},x_{2},...,x_L)=\sum_{r_1,...,r_L\geq
0}\mathcal{H}_{r_1,...,r_L}(\mathcal{L}^{\vec{\tau}};q,a)x_1^{r_1}\cdots x_L^{r_{L}}.
\end{align}
Then the Ooguri-Vafa's conjecture \cite{OV} for knot can be generalized to the case of framed link as follows. 
\begin{conjecture} \label{conjecture-framedlink}
The generating function (\ref{formula-Zlink}) can be rewritten in the following form 
\begin{align*}  
    \mathcal{Z}(\mathcal{L}^{\vec{\tau}};q,a;x_{1},x_{2},...,x_L)=\exp\left(\sum_{d\geq 1}\sum_{\substack{r_1,r_2,...,r_L\geq 0\\ (r_1,..,r_L)\neq (0,..,0)}}\frac{1}{d}\mathfrak{f}_{r_1,...,r_L}(\mathcal{L}^{\vec{\tau}};a^d,q^d)x_1^{dr_1}\cdots x_L^{dr_L})\right). 
\end{align*}
Suppose $r_1,r_2,...,r_L\geq 0$ and $(r_1,..,r_L)\neq (0,..,0)$, let $k=\#\{r_i|r_i\neq 0, i=1,...,L\}$, and hence $1\leq k\leq L$. There exists integers $N_{(r_1,...,r_L),i,j}(\mathcal{L}^{\vec{\tau}})$, that vanish for large $|i|,|j|$, such that  
\begin{align}
    \mathfrak{f}_{r_1,...,r_L}(\mathcal{L}^{\vec{\tau}};q,a)=(q^{\frac{1}{2}}-q^{-\frac{1}{2}})^{k-2}\sum_{i,j\in \mathbb{Z}}N_{(r_1,...,r_L),i,j}(\mathcal{L}^{\vec{\tau}})a^{\frac{i}{2}}q^{\frac{j}{2}}.
\end{align}
\end{conjecture}

Based on the Conjecture \ref{conjecture-framedlink}, we introduce the polynomial of $a^{\pm \frac{1}{2}}$ and $q^{\pm \frac{1}{2}}$ as follows
\begin{align}
p_{r_1,...r_L}(\mathcal{L}^{\vec{\tau}};q,a)=\sum_{i,j\in \mathbb{Z}}N_{(r_1,...,r_L),i,j}(\mathcal{L}^{\vec{\tau}})a^{\frac{i}{2}}q^{\frac{j}{2}}.
\end{align}
By the strongly integral property of colored HOMFLYPT invariants proved in \cite{Zhu23},  we obtain
\begin{align}
p_{r_1,...r_L}(\mathcal{L}^{\vec{\tau}};q,a)=\sum_{i,j\in \mathbb{Z}}N_{(r_1,...,r_L),i,j}(\mathcal{L}^{\vec{\tau}})a^{\frac{i}{2}}q^{\frac{j}{2}}\in
a^{\frac{1}{2}\epsilon_1}q^{\frac{1}{2}\epsilon_2} \mathbb{Z}[a^{\pm 1},q^{\pm 1}]
\end{align}
where $\epsilon_1,\epsilon_2\in \{0,1\}$ is determined by $$ |r_1+\cdots +r_L|\equiv
\epsilon_1 \mod 2 $$ and $$ |r_1+\cdots +r_L|+k\equiv
\epsilon_2 \mod 2. $$

Finally, we compute the Ooguri-Vafa invariants and BPS invariants of framed Whitehead links and Borromean rings by using the explicit formulas of the colored HOMFLYPT invariants for them conjectured in \cite{GNSSS15} and proved in \cite{CLZ}. These calculations confirm Conjecture \ref{conjecture-framedlink}.

The rest of this article is organized as follows. In Section \ref{section-BPSfromframedknots}, we review the background for the large $N$ duality of Chern-Simons and topological string theories and introduce the BPS invariants for the framed knot.  The relations between algebraic curves and BPS invariants of the framed knots are then presented in Section \ref{section-dualAextremal}. We give the framing change formula for the dual $A$-polynomial, as well as the extremal $A$-polynomial of a framed knot. In Section \ref{section-example}, we illustrate the detailed computations for the extremal $A$-polynomial and extremal BPS invariant of framed unknot $U^\tau$ and framed twist knots $\mathcal{K}_p^{\tau}$. Section \ref{section-link} is devoted to the cases of framed links.

\section{BPS invariants of framed knots} \label{section-BPSfromframedknots}
Let us consider the open sector of topological A-model of a
Calabi-Yau 3-fold $X$ with a Lagrangian submanifold $\mathcal{D}$
of dim $H_{1}(\mathcal{D},\mathbb{Z})=1$. The open sector
topological A-model can be described by holomorphic maps $\phi$ from
open Riemann surface of genus $g$ and $l$-holes $\Sigma_{g,l}$ to
$X$, with Dirichlet condition specified by $\mathcal{D}$.  These
holomorphic maps are referred as open string instantons. More
precisely, an open string instanton is a holomorphic map $\phi:
\Sigma_{g,l}\rightarrow X$ such that $\partial \Sigma_{g,l}=\cup
_{i=1}^{l}\mathcal{C}_i\rightarrow \mathcal{D}\subset X$ where the
boundary $\partial \Sigma_{g,l}$ of $\Sigma_{g,l}$ consists of $l$
connected components $\mathcal{C}_i$ mapped to Lagrangian
submanifold $\mathcal{D}$ of $X$. Therefore, the open string
instanton $\phi$ is described by the following two different kinds
of data: the first is the ``bulk part" which is given  by
$\phi_*[\Sigma_{g,l}]=Q\in H_2(X,\mathcal{D})$, and the second is
the ``boundary part" which is given by
$\phi_*[\mathcal{C}_i]=w_i\gamma$, for $i=1,..l$,
where $\gamma$ a generator of
$H_{1}(\mathcal{D},\mathbb{Z})$ and $w_i\in \mathbb{Z}$. Let
$w=(w_1,..,w_l)$, we expect that there exists the corresponding open
Gromov-Witten invariant $K_{w,g,Q}$ interpreted as the virtual counting of the open string instantons determined by the data
 $w, Q$ in the genus $g$. See \cite{LS,KL} for mathematical
aspects of defining these invariants in special cases.

We take all $w_i\geq 1$ as in \cite{MV}, and we use the notation of
partitions and symmetric functions \cite{Mac}. We denote by
$\mathcal{P}$ the set of all partitions, including the empty
partition $0$, and by $\mathcal{P}_+$ the set of nonzero partitions.
Let $\mathbf{x}=\{x_1,x_2,...\}$ be the set of infinitely many
independent variables. For $n\geq 0$, let
$p_n(\mathbf{x})=\sum_{i\geq 1}x_i^{n}$ be a symmetric power sum
function. For a partition $\mu\in \mathcal{P}_{+}$, set
$p_{\mu}(\mathbf{x})=\prod_{i=1}^{h}p_{\mu_i}(\mathbf{x})$. 

 The total
free energy and partition function of open topological string on $(X,\mathcal{D})$
are defined as follows
\begin{align} \label{partion-freeenergy}
F^{(X,\mathcal{D})}(\mathbf{x};g_s,\omega)&=-\sum_{g\geq
0}\sum_{\mu \in \mathcal{P}\setminus
\{0\}}\frac{\sqrt{-1}^{l(\mu)}}{|Aut(\mu)|}g_{s}^{2g-2+l(\vec{\mu})}
\sum_{Q\neq 0}K_{\mu,g,Q}^{(X,\mathcal{D})}a^{\frac{Q}{2}} p_{\mu}(\mathbf{x}),\\\nonumber
Z^{(X,\mathcal{D})}(\mathbf{x};g_s,a)&=\exp\left(F^{(X,\mathcal{D})}(\mathbf{x};g_s,a)\right).
\end{align}

\subsection{Large $N$ duality for framed knots}
In his seminal paper \cite{W1}, E. Witten introduced a new
topological invariant of a 3-manifold $M$  as a partition function
of quantum Chern-Simons theory. Let $G$ be a compact gauge group
that is a Lie group, and $M$ be an oriented three-dimensional
manifold. Let $\mathcal{A}$ be a $\mathfrak{g}$-valued connection on
$M$ where $\mathfrak{g}$ is the Lie algebra of $G$. The Chern-Simons
\cite{CS} action is given by
\begin{align*}
S(\mathcal{A})=\frac{k}{4\pi}\int_{M}Tr\left(\mathcal{A}\wedge
d\mathcal{A}+\frac{2}{3}\mathcal{A}\wedge\mathcal{A}
\wedge\mathcal{A}\right)
\end{align*}
where $k$ is an integer called the level.

Chern-Simons partition function is defined as the path integral in
quantum field theory
\begin{align*}
Z^G(M;k)=\int e^{i S(A)}D \mathcal{A}
\end{align*}
where the integral is over the space of all $\mathfrak{g}$-valued
connections $\mathcal{A}$ on $M$. Although it is not rigorous,
Witten \cite{W1} developed some techniques to calculate such
invariants.

If the three-manifold $M$ contains a link $\mathcal{L}$, we let
$\mathcal{L}$ be an $L$-component link denoted by
$\mathcal{L}=\bigsqcup _{j=1}^L\mathcal{K}_j$. Define
$$W_{R_j}(\mathcal{K}_j)=Tr_{R_j}Hol_{\mathcal{K}_j}(\mathcal{A})$$ which is the trace of holomony
along $\mathcal{K}_j$ taken in representation $R_j$. Then Witten's
invariant of the pair $(M,\mathcal{L})$ is given by
\begin{align*}
Z^{G}(M,\mathcal{L};\{R_j\};k)=\int e^{iS(\mathcal{A})}\prod_{j=1}^L
W_{R_j}(\mathcal{K}_j)D\mathcal{A}.
\end{align*}

Later, Reshetikhin
and Turaev \cite{RT1,RT2} developed a systematic way to construct
the above invariants by using the representation theory of quantum
groups. Their construction led to the definition of colored
HOMFLYPT invariants \cite{LZ}, which can be viewed as the
large $N$ limit of the quantum $U_{q}(sl_N)$ invariants. Usually, we
use the notation $H_{\lambda^1,..,\lambda^L}(\mathcal{L};q,a)$ to
denote the (framing-independent) colored HOMFLY-PT invariants for a
(oriented) link $\mathcal{L}=\bigsqcup _{j=1}^L\mathcal{K}_j$, where
each component $\mathcal{K}_j$ is colored by an irreducible
representation $V_{\lambda^j}$ of $U_{q}(sl_N)$ for $N$ large enough. Some basic
structures for $H_{\lambda^1,..,\lambda^L}(\mathcal{L};q,a)$ were
proved in \cite{GLL18,LP1,LP2,Zhu2}. It is difficult to obtain an explicit
formula for a given link with irreducible representations. We refer to \cite{LZ} for an explicit formula for torus
links. In a recent work \cite{CLZ} joined with Q. Chen and K. Liu, the second author obtained  explicit formulas for the colored HOMFLYPT invariants (under symmetric representations)
of double twist knots, twisted Whitehead links and Borromean rings, which were conjectured in \cite{MMM,NRZ,GNSSS15} previously.  

In another fundamental work of Witten \cite{W2}, the $SU(N)$
Chern-Simons gauge theory on a
three-manifold $M$ was interpreted as an open topological string theory on $%
T^*M$ with $N$ topological branes wrapping $M$ inside $T^*M$.
Furthermore, Gopakumar and Vafa \cite{GV2} conjectured that the
large $N$ limit of $SU(N)$ Chern-Simons gauge theory on $S^3$ is
equivalent to closed topological string theory on the resolved
conifold. Furthermore, Ooguri and Vafa \cite{OV} generalized the
above construction to the case of a knot $\mathcal{K}$ in $S^3$.
They introduced the Chern-Simons partition function
$Z_{CS}^{(S^3,\mathcal{K})}(q,a,\mathbf{x})$ for $(S^3,\mathcal{K})$
which is a generating function of the colored HOMFLY-PT invariants
in all irreducible representations.
\begin{align} \label{chernsimonspartition}
Z_{CS}^{(S^3,\mathcal{K})}(q,a,\mathbf{x})=\sum_{\lambda\in
\mathcal{P}}H_{\lambda}(\mathcal{K},q,a)s_{\lambda}(\mathbf{x}),
\end{align}
where $s_{\lambda}(\mathbf{x})$ is the Schur function \cite{Mac}.

Ooguri and Vafa \cite{OV} conjectured that for any knot
$\mathcal{K}$ in $S^3$, there exists a corresponding Lagrangian
submanifold $\mathcal{D}_{\mathcal{K}}$, such that the Chern-Simons
partition function
 is equal to the open topological string
partition function
 on
$(X,\mathcal{D}_{\mathcal{K}})$. They have established this duality
for the case of a trivial knot $U$ in $S^3$, and the link case was
further discussed in \cite{LMV}.

 In general,  we first
must find a way to construct the Lagrangian submanifold
$\mathcal{D}_\mathcal{L}$ corresponding to the link $\mathcal{L}$ in
geometry. See \cite{LMV,Koshkin,Tau,DSV} for the constructions for
some special links. Furthermore, if the Lagrangian submanifold
$\mathcal{D}_{\mathcal{L}}$ is constructed, then we need to compute
the open sting partition function under this geometry.  Aganagic and Vafa \cite{AV} introduced the
special Lagrangian submanifold in toric Calabi-Yau 3-fold which we
call Aganagic-Vafa A-brane (AV-brane) and studied its mirror
geometry, then they computed the counting of the holomorphic disc end on
AV-brane by using the idea of mirror symmetry. Moreover, Aganagic
and Vafa surprisingly found the computation by using mirror symmetry
and the result from the Chern-Simons knot invariants \cite{OV} are
matched. Furthermore, in \cite{AKV}, Aganagic, Klemm and Vafa
investigated the integer ambiguity appearing in the disc counting
and discovered that the corresponding ambiguity in Chern-Simons
theory was described by the framing of the knot. They checked that
the two ambiguities match for the case of the unknot, by comparing
the disk amplitudes on both sides.

Then, Mari\~no and Vafa \cite{MV} generalized the large $N$ duality
to the case of knots with arbitrary framing. More precisely,
for a framed knot $\mathcal{K}^\tau$ with framing $\tau\in
\mathbb{Z}$, we define the framed colored HOMFLYPT invariants
$\mathcal{K}^\tau$ as follows
\begin{align} \label{framedknotformula}
\mathcal{H}_\lambda(\mathcal{K}^\tau,q,a)=(-1)^{|\lambda|\tau}q^{\frac{\kappa_\lambda\tau}{2}}H_{\lambda}(\mathcal{K},q,a),
\end{align}
where
$\kappa_\lambda=\sum_{i=1}^{l(\lambda)}\lambda_i(\lambda_i-2i+1)$.

 The framed Chern-Simon partition function for
$(S^3,\mathcal{K}^\tau)$ is given by
\begin{align} \label{partitionfunctionframedunknot}
\mathcal{Z}_{CS}^{(S^3,\mathcal{K}^\tau)}(q,a;\mathbf{x})=\sum_{\lambda\in
\mathcal{P}}\mathcal{H}_\lambda(\mathcal{K}^\tau,q,a)s_{\lambda}(\mathbf{x}).
\end{align}
We let $\hat{X}:=\mathcal{O}(-1)\oplus\mathcal{O}(-1)\rightarrow
\mathbb{P}^1$ be the resolved conifold and $D_{\mathcal{K}^\tau}$ be the
corresponding Lagrangian submanifold. The open string partition function for
$(\hat{X},D_{\mathcal{K}^\tau})$ has the structure
\begin{align} \label{formula-STpartitonfunction}
Z_{str}^{(\hat{X},D_{\mathcal{K}^\tau})}(g_s,a;\mathbf{x})=\exp\left(
-\sum_{g\geq 0,\mu\in \mathcal{P}\setminus\{0\}}\frac{\sqrt{-1}^{l(\mu)}}{|Aut(\mu)|}g_s^{2g-2+l(\mu)}F_{\mu,g}^{
\tau}(a)p_{\mu}(\mathbf{x})\right)
\end{align}
where $ F^{\tau}_{\mu,g}(a)=\sum_{Q\in
\mathbb{Z}}K^{\tau}_{\mu,g,Q}a^\frac{Q}{2} $ and $K_{\mu,g,Q}^{\tau}$ is the expected
open Gromov-Witten invariants for $(\hat{X},D_{\mathcal{K}^\tau})$.

Therefore, the large $N$ duality for framed knot $\mathcal{K}^\tau$ is given by the
following identity:
\begin{align} \label{formula-CSTS}
\mathcal{Z}_{CS}^{(S^3,\mathcal{K}^\tau)}(q,a;\mathbf{x})=Z_{str}^{(\hat{X},D_{\mathcal{K}^\tau})}(g_s,a;\mathbf{x})
\end{align}
where $q=e^{ig_s}$.

\subsection{BPS invariants for framed knots}
Based on the large $N$ duality, Ooguri-Vafa \cite{OV} conjectured that for any knot $\mathcal{K}$, the generating function $Z_{CS}^{(S^3,\mathcal{K})}(q,a;\mathbf{x})$ has the following form
\begin{align}
    Z_{CS}^{(S^3,\mathcal{K})}(q,a;\mathbf{x})=\exp\left(\sum_{d\geq 1}\sum_{\mu\in \mathcal{P}^+}\frac{1}{d}f_{\mu}(\mathcal{K};a^d,q^d)s_\mu(\mathbf{x}^d)\right),
\end{align}
where the functions $f_{\mu}(\mathcal{K};a,q)$ has the form
\begin{align}
  f_{\mu}(\mathcal{K};a,q)=\sum_{i,j}\frac{N_{\mu,i,j}(\mathcal{K})a^{\frac{i}{2}}q^{\frac{j}{2}}}{q^{\frac{1}{2}}-q^{-\frac{1}{2}}},   
\end{align}
where $N_{\mu,i,j}(\mathcal{K})$ are the famous $BPS$ degeneracies, or Ooguri-Vafa invariants,\footnote{In some literatures, such as \cite{GKS}, refer $N_{\mu,i,j}(\mathcal{K})$ as LMOV (Labastida-Mari\~no-Ooguri-Vafa) invariants. In fact, Ooguri and Vafa \cite{OV} first proposed BPS degeneracies in the form $N_{\mu,i,j}(\mathcal{K})$. Later, Labastida, Mari\~no and Vafa \cite{LMV} proposed the refined form of BPS degeneracies $N_{\mu,g,Q}(\mathcal{K})$, which are referred to as LMOV invariants \cite{LP1}} in particular, they are conjectured to be integer. Furthermore, for a fixed $\mu$, there is a finite range of $i$ and $j$ for which $N_{\mu,i,j}(\mathcal{K})$ are non-zero.

We stress that in the definition of the colored HOMFLYPT invariants of a knot, we choose as our variables $q^{\frac{1}{2}},a^{\frac{1}{2}}$ rather than $q,a$ as in \cite{GKS}, in order to keep the consistency with the original physical literature \cite{OV,LMV}.

Now, let us consider a framed knot $\mathcal{K}^\tau$, based on Marino and Vafa's work , it is also conjectured that the framed Chern-Simons partition    $\mathcal{Z}_{CS}^{(S^3,\mathcal{K}^\tau)}(q,a;\mathbf{x})$ has the following form
\begin{align} \label{formula-CSpartitionfunction}
    \mathcal{Z}_{CS}^{(S^3,\mathcal{K}^\tau)}(q,a;\mathbf{x})=\exp\left(\sum_{d\geq 1}\sum_{\mu\in \mathcal{P}^+}\frac{1}{d}\mathfrak{f}_{\mu}(\mathcal{K}^\tau;a^d,q^d)s_\mu(\mathbf{x}^d)\right),
\end{align}
where the functions $\mathfrak{f}_{\mu}(\mathcal{K}^\tau;a,q)$ are given by
\begin{align} \label{formula-ftau}
  \mathfrak{f}_{\mu}(\mathcal{K}^\tau;a,q)=\sum_{i,j}\frac{N_{\mu,i,j}(\mathcal{K}^\tau) a^{\frac{i}{2}}q^{\frac{j}{2}}}{q^{\frac{1}{2}}-q^{-\frac{1}{2}}},
\end{align}
where $N_{\mu,i,j}(\mathcal{K}^\tau)$ are integers, and named as framed Ooguri-Vafa invariants.  

In what follows we are interested in the case of one-dimensional source, i.e.  $\mathbf{x}=(x,0,...)$.  Then $s_{\lambda}(\mathbf{x})\neq 0$ only when $\lambda=(r)$ which is a single row with $n$ boxes, and $s_{(r)}(\mathbf{x})=x^r$. We let  
\begin{align} \label{formula-framedH}
\mathcal{H}_{r}(\mathcal{K}^\tau,a,q):=\mathcal{H}_{(r)}(\mathcal{K}^\tau,a,q)=(-1)^\tau q^{\frac{r(r-1)\tau}{2}}H_r(\mathcal{K};a,q).
\end{align}

In this setting, formula (\ref{formula-CSpartitionfunction}) is reduced to the following
\begin{align}  \label{formula-reducedCSknot}
     \mathcal{Z}_{CS}^{(S^3,\mathcal{K}^\tau)}(x,q,a)=\sum_{r\geq 0}\mathcal{H}_{r}(\mathcal{K}^\tau,q,a)x^r=\exp\left(\sum_{d\geq 1}\sum_{r\geq 1}\frac{1}{d}\mathfrak{f}_{r}(\mathcal{K}^\tau;a^d,q^d)x^{rd})\right).
\end{align}
By M\"obius inversion formula, we obtain 
\begin{align} \label{formula-fn}
\mathfrak{f}_r(\mathcal{K}^{\tau};a,q)=\sum_{d|r}\frac{\mu(d)}{d}\sum_{|\mu|=r/d}\frac{(-1)^{l(\mu)}(l(\mu)-1)!}{|Aut(\mu)|}\prod_{i=1}^{l(\mu)}\mathcal{H}_{\mu_i}(\mathcal{K}^\tau;a^d,q^d),    
\end{align}
where the notations $l(\mu)$ and $Aut(\mu)$ denote the length and automorphism group of the partition $\mu$.  From the expression (\ref{formula-fn}), it follows that $\mathfrak{f}_r(\mathcal{K}^\tau;a,q)\in \mathbb{Q}(a,q)$.  Ooguri-Vafa's conjecture \cite{OV} extended to framed knot $\mathcal{K}_\tau$ asserts that $\mathfrak{f}_r(\mathcal{K}^{\tau};a,q)$ can be expressed as a finite sum
\begin{align} \label{formula-fr-integral}
\mathfrak{f}_{r}(\mathcal{K}^\tau;a,q)=\sum_{i,j}\frac{N_{r,i,j}(\mathcal{K}^\tau) a^{\frac{i}{2}}q^{\frac{j}{2}}}{q^{\frac{1}{2}}-q^{-\frac{1}{2}}}, \  N_{r,i,j}(\mathcal{K}^\tau) \in \mathbb{Z}. 
\end{align}
A natural question is how to extract the Ooguri-Vafa invariant $N_{r,i,j}(\mathcal{K}^\tau)$ from the generating function (\ref{formula-CSpartitionfunction}). As in \cite{GKS}, we introduce 
\begin{align}
b_{r,i}(\mathcal{K}^\tau)=\sum_{j}N_{r,i,j}(\mathcal{K}^\tau)    
\end{align}
which can be interpreted as the classical BPS degeneracies.  

\subsection{String theory interpretation of the BPS invariants $b_{r,i}(\mathcal{K}^\tau)$}
In open string theory, these invariants $b_{r,i}(\mathcal{K}^\tau)$ determine the holomorphic disk amplitudes of $D$-branes conjecturally associated with the framed knot $\mathcal{K}^\tau$. In fact, taking $\mathbf{x}=(x,...)$ in formula (\ref{formula-STpartitonfunction}) and (\ref{formula-CSTS}), we obtain 
\begin{align}
 \mathcal{Z}^{(S^3,\mathcal{K}^\tau)}(q,a;x)=\exp\left(-\sum_{g\geq 0,n\geq 1}\sqrt{-1}\sum_{Q}K_{n,g,Q}^\tau g_s^{2g-1} a^{\frac{Q}{2}}x^n \right)   
\end{align}

Comparing with formula (\ref{formula-CSpartitionfunction}), for any $r\geq 1$,  we have
\begin{align} 
 \sum_{d|r}\frac{1}{d}\mathfrak{f}_{r/d}(\mathcal{K}^\tau;a^d,q^d)= -\sum_{g\geq 0}\sqrt{-1}\sum_{Q}K_{r,g,Q}^\tau g_s^{2g-1} a^{\frac{Q}{2}}.   
\end{align}
Substituting the expression (\ref{formula-fr-integral}) to above formula, we obtain 
\begin{align}
\sum_{d|r,d|i}\frac{1}{d}\sum_{j}\frac{N_{\frac{r}{d},\frac{i}{d},j}(\mathcal{K}^\tau)q^{\frac{dj}{2}}}{q^{\frac{d}{2}}-q^{-\frac{d}{2}}}=-\sqrt{-1}\sum_{g\geq  0,i}K_{r,g,i}^{\tau}g_{s}^{2g-1}.    
\end{align}
Multiply it by $g_s$ and taking the limit $g_s\rightarrow 0$, we get
\begin{align} \label{formula-GromovWitten}
K_{r,0,i}^\tau=\sum_{d|r,d|i}\sum_{j}\frac{N_{\frac{r}{d},\frac{i}{d},j}(\mathcal{K}^\tau)}{d^2}=\sum_{d|r,d|i}\frac{1}{d^2}b_{\frac{r}{d},\frac{i}{d}}(\mathcal{K}^\tau),
\end{align}
since $
 \lim_{g_s\rightarrow 0}\frac{g_sq^{\frac{dj}{2}}}{q^{\frac{d}{2}}-q^{-\frac{d}{2}}} =\frac{1}{\sqrt{-1}d}.  
$

Notice that the geometric meaning of $K_{r,0,i}^\tau$ is the counting of the holomorphic disk invariants by the description of the open Gromov-Witten given at the beginning of Section \ref{section-BPSfromframedknots}. Hence, formula (\ref{formula-GromovWitten}) shows that $K_{n,0,i}^\tau$ is determined by BPS invariants $b_{r,i}(\mathcal{K}^\tau)$. On the other hand, by M\"obius inversion formula, we also have
\begin{align}  \label{formula-bps-GW}
    b_{r,i}(\mathcal{K}^\tau)=\sum_{d|r,d|i}\frac{\mu(d)}{d^2}K_{\frac{r}{d},0,\frac{i}{d}}^{\tau}. 
\end{align}

\section{Dual $A$-polynomial and extremal invariants for framed knots} \label{section-dualAextremal}
 
A natural question is how to compute the BPS invariant $b_{r,i}(\mathcal{K}^\tau)$. We describe a method introduced in \cite{GKS} to compute them using an algebraic curve, namely dual $A$-polynomial, associated with a knot. Let us briefly review it following the states in \cite{GKS}. Recall that an one-parameter sequence $\{f_n|\in \mathbb{N}\}$ of polynomials in $\mathbb{Z}[a^{\pm}](q^{\frac{1}{2}})$ is $q$-holonomic if there exists $d\in \mathbb{N}$ and polynomials $a_l(u,v)\in \mathbb{Z}[u,v]$ for $0\leq l\leq d$ such that 
\begin{align}
\sum_{l=0}^da_{l}(q^{\frac{1}{2}},q^{\frac{n}{2}})f_{n+l}=0   \ \text{for all $n\geq  0$}.   
\end{align}
Note that this recurrence can be encoded in a polynomial $\hat{A}(M,L,a,q)$ with $LM=q^{\frac{1}{2}}ML$, where $M$ and $L$ are considered as operators that act on the sequence $\{f_n\}$ by 
\begin{align}
 M f_n=q^\frac{n}{2}f_n, \ Lf_n=f_{n+1}   
\end{align}

It is proved in \cite{GLL18} that the HOMFLYPT polynomial of a link colored by partitions with a fixed number of rows is a $q$-holonomic function. For the framed knot $\mathcal{K}^\tau$, the framed colored HOMFLYPT invariant $\mathcal{H}_{n}(\mathcal{K}^\tau;a,q)$ satisfies a linear $q$-difference equation of the form
\begin{align}
\hat{A}(M,L,a,q)\mathcal{H}_{r}(\mathcal{K}^\tau,a,q)=0.     
\end{align}

The corresponding 3-variable polynomial $A_{\mathcal{K}^\tau}(M,L,a)=A(M,L,a,1)$ defines a family of algebraic curves parameterized by $a$
which is referred to as the $a$-deformed $A$-polynomial of $\mathcal{K}^\tau$ in \cite{AV12} (note that the variable $a$ coincides with $Q$ in \cite{AV12}).  The dual $A$-polynomial of $A_{\mathcal{K}^\tau}(M,L,a)$ is defined by \cite{GKS}
\begin{align} \label{formula-AdualA}
   \mathcal{A}_{\mathcal{K}^\tau}(x,y,a)=A_{\mathcal{K}^\tau}(y,x^{-1},a). 
\end{align}

Using the conjectural integrality formula (\ref{formula-ftau}), Proposition 1.1 in \cite{GKS} gives the following formula for the computation of the BPS invariant $b_{r,i}(\mathcal{K}^\tau)$. 

Given a framed knot $\mathcal{K}^\tau$, 
suppose that $y=y(x,a;\tau)$ is a solution of the dual $A$-polynomial $\mathcal{A}_{\mathcal{K}^\tau}(x,y,a)=0$, then the BPS invariants $b_{r,i}(\mathcal{K}^\tau)$ can be read from the following formula
\begin{align} \label{formula-bpsfromcurve}
  x\frac{\partial_x y}{y}=-\frac{1}{2}\sum_{r,i}r^2b_{r,i}(\mathcal{K}^\tau)\frac{x^ra^{\frac{i}{2}}}{1-x^ra^{\frac{i}{2}}}.  
\end{align}
Hence,  if we have
\begin{align}
    x\frac{\partial_x y}{y}=c_{r,i}(\mathcal{K}^\tau)x^ra^{\frac{i}{2}},
\end{align}
 then, by M\"obius inversion formula, we obtain
\begin{align} \label{formula-bri}
 b_{r,i}(\mathcal{K}^\tau)=\frac{2}{r^2}\sum_{d|r,i}\mu(d)c_{\frac{r}{d},\frac{i}{d}}(\mathcal{K}^\tau).   
\end{align}

\subsection{Computations of the dual $A$-polynomial $\mathcal{A}_{\mathcal{K}^\tau}(x,y,a)$}

Suppose $A_{\mathcal{K}}(x,y,a)$ is the $a$-deformed $A$-polynomial of $\mathcal{K}$. By its definition, $A_{\mathcal{K}}(x,y,a)$ is determined by the classical limit of recursions for colored HOMFLYPT invariants $H_n(\mathcal{K};q,a)$.  

Assume that $H_{n}(\mathcal{K};q,a)$ has the form of multiple summation $\sum_{0\leq k_1,...,k_l\leq n}\cdots$ and suppose $q=e^{h}$, we introduce the variables $\tilde{x}=q^n$ and $z_i=q^{k_i}$ for $i=1,..,l$.  Then $n=\frac{\log(x)}{h}$ and $k_i=\frac{\log(z_i)}{h}$ for $i=1,...,l$. By using the following asymptotic expansion formula for quantum dilogarithm function 
\begin{align}
    (z;q)_k\sim_{h\rightarrow 0}e^{\frac{1}{h}(\text{Li}_2(z)-\text{Li}_2(zq^{k-1}))}, 
\end{align}
we obtain the potential function $W_{\mathcal{K}}(\tilde{x},z_1,...,z_l;a)$ for $\mathcal{K}$. The critical point equations for  $W_{\mathcal{K}}(\tilde{x},z_1,...,z_l;a)$ is given by 
\begin{align}
 \exp\left(z_i\frac{\partial W_{\mathcal{K}}(\tilde{x},z_1,...,z_l;a)}{\partial z_i}\right)=1,  \ \ i=1,...,l.   
\end{align}

Then the $a$-deformed $A$-polynomial  $A_{\mathcal{K}}(x,y,a)=0$ for $\mathcal{K}$ is obtained by eliminating $z_1^0,...,z_l^0$ from 
\begin{align} \label{formula-yx}
 y=\exp\left(\tilde{x}\frac{\partial W_{\mathcal{K}}(\tilde{x},z_1^0,...,z_l^0;a)}{\partial \tilde{x}}\right)   
\end{align}
using the saddle point equations, and substituting $\tilde{x}=x^2$ in the resulting equation.

By using formula (\ref{formula-framedH}), it is easy see that 
\begin{align}
    W_{\mathcal{K}^\tau}(\tilde{x},z_1,...,z_l;a)=\tau \pi\sqrt{-1}\log \tilde{x}+\tau\frac{(\log \tilde{x})^2}{2}+W_{\mathcal{K}}(\tilde{x},z_1,...,z_l;a).
\end{align}
Hence, the formula (\ref{formula-yx}) gives  
$A_{\mathcal{K}^\tau}(x,y,a)=A_{\mathcal{K}}(x,(-1)^\tau x^{-2\tau}y,a)$. Therefore, 
\begin{align}
    \mathcal{A}_{\mathcal{K}^\tau}(x,y,a)=A_{\mathcal{K}^\tau}(y,x^{-1},a)=A_{\mathcal{K}}(y,(-1)^\tau y^{-2\tau} x^{-1},a). 
\end{align}

Alternatively, we can derive the dual $A$-polynomial for framed knot $\mathcal{K}^\tau$ from the dual $A$-polynomial of knot $\mathcal{K}^0$ as follows 
\begin{align} \label{formula-AtauandA}
\mathcal{A}_{\mathcal{K}^\tau}(x,y,a)=\mathcal{A}_{\mathcal{K}^0}((-1)^\tau y^{2\tau} x,y,a), 
\end{align}
where $\mathcal{A}_{\mathcal{K}^0}(x,y,a)$ is coincided with the dual $A$-polynomial $\mathcal{A}_{\mathcal{K}}(x,y,a)$ introduced in \cite{GKS}.

\subsection{Extremal invariants and extremal $A$-polynomial}
In \cite{GKS}, Garoufalidis et al introduced the extremal BPS invariants of a knot. It's definition also works for a framed knot $\mathcal{K}_\tau$.  The basic setting is that the minimal and maximal exponent of $\mathcal{H}_r(\mathcal{K}^\tau,a,q)$ with respect to $a$ is a quasi-linear function of $r$ for large enough $r$, since $\mathcal{H}_r(\mathcal{K}^\tau,a,q)$ is a $q$-holonomic sequence.  For simplicity, we only consider the special case that the framed colored HOMFLYPT invariants $\mathcal{H}_r(\mathcal{K}^\tau,a,q)$  satisfy  

\begin{align}
\mathcal{H}_r(\mathcal{K}^\tau;a,q)=\sum_{i=rc_-}^{rc_+}a^{\frac{i}{2}}h_{r,i}(\mathcal{K}^\tau,q)    
\end{align}
for some integers $c_\pm$ and for every natural number $r$, where $h_{r,r\cdot c_\pm}(\mathcal{K}^\tau,q)\neq 0$. This is a large class of knots, in particular two-bridge knots and torus knots have this property.  Using this property, one can define the extremal BPS invariants 
\begin{align}
b_{r}^\pm(\mathcal{K}^\tau)=b_{r,r\cdot c_\pm}(\mathcal{K}^\tau)=\sum_{j}N_{r,r\cdot c_\pm,j}(\mathcal{K}^\tau).    
\end{align}
Moreover, the extremal analogues of the dual $A$-polynomial (\ref{formula-AtauandA}) is  defined as distinguished, irreducible factors in 
\begin{align}
    \mathcal{A}_{\mathcal{K}^\tau}(xa^{\pm c_\pm},y,a)|_{a^{\mp 1}\rightarrow 0}
\end{align}
that determine the extremal BPS invariant $b_{r,i}(\mathcal{K}^\tau)$. These curves are called the extremal $A$-polynomials and are denoted by $\mathcal{A}^\pm_{\mathcal{K}^\tau}(x,y)$. Similar to the formula  (\ref{formula-AtauandA}), we have 
\begin{align} \label{formula-AdualApm}
  \mathcal{A}_{\mathcal{K}^\tau}^{\pm }(x,y)=\mathcal{A}_{\mathcal{K}^0}^\pm((-1)^\tau y^{2\tau} x,y).   
\end{align}

\section{Examples}  \label{section-example}
In this section, we illustrate the claims and ideas presented previously in several examples.

\subsection{Framed unknot}
First, let us illustrate the computations for the framed unknot $U^\tau$. We
first compute the non-commutative a-deformed $A$-polynomial (it is
called the Q-deformed A-polynomial in \cite{AV12}, the variable $Q$
in \cite{AV12} is the variable $a$ here) for $U^\tau$.

Recall that the colored HOMFLYPT invariants
$H_{n}(U;q,a)$ of unknot $U$ is given by
\begin{align} \label{formula-Hu}
H_{n}(U;q,a)=\frac{a^{\frac{1}{2}}-a^{-\frac{1}{2}}}{q^{\frac{1}{2}}-q^{-\frac{1}{2}}}\cdots
\frac{a^{\frac{1}{2}}q^{\frac{n-1}{2}}-a^{-\frac{1}{2}}q^{\frac{-n-1}{2}}}{q^{\frac{n}{2}}-q^{-\frac{n}{2}}}.
\end{align}

By formula (\ref{framedknotformula}), the framed colored HOMFLYPT
invariants for the framed unknot $U^\tau$ with framing $\tau\in \mathbb{Z}$
is given by
\begin{align}
\mathcal{H}_n(U^\tau;q,a)=(-1)^{n\tau}q^{\frac{n(n-1)}{2}\tau}H_n(U;q,a).
\end{align}
Then we obtain the recursion for $\mathcal{H}_n(U^\tau;q,a)$ as follows
\begin{align} \label{frameunknotrecursion}
(-1)^\tau(q^{n+1}-1)\mathcal{H}_{n+1}(U^\tau;q,a)-(a^{\frac{1}{2}}q^{n+\frac{1}{2}}-a^{-\frac{1}{2}}q^{\frac{1}{2}})q^{n\tau}\mathcal{H}_n(U^\tau;q,a)=0.
\end{align}
It implies that the noncommutative $a$-deformed $A$-polynomial for $U_\tau$ is given by
\begin{align*}
\hat{A}_{U^\tau}(M,L,q;a)=(-1)^\tau(qM^2-1)L-M^{2\tau}(a^{\frac{1}{2}}q^{\frac{1}{2}}M^2-a^{-\frac{1}{2}}q^{\frac{1}{2}}).
\end{align*}
Then we obtain the a-deformed $A$-polynomial of $U^\tau$:
\begin{align*}
A_{U^\tau}(M,L;a)=\lim_{q\rightarrow
1}\hat{A}(M,L,q;a)=(-1)^\tau(M^2-1)L-M^{2\tau}(a^{\frac{1}{2}}M^2-a^{-\frac{1}{2}}).
\end{align*}
By formula (\ref{formula-AdualA}), the dual $a$-deformed $A$-polynomial of $U^\tau$ is given by
\begin{align}
\mathcal{A}_{U^\tau}(x,y,a)=A_{U^\tau}(y,x^{-1};a).
\end{align}
Hence, we have
\begin{align} \label{formula-AUtau}
\mathcal{A}_{U^\tau}(x,y;a)=y^2-1-a^{-\frac{1}{2}}(-1)^\tau
xy^{2\tau}(ay^2-1).
\end{align}

By formula (\ref{formula-Hu}), we obtain $c_\pm=\pm 1$, the extremal $A$-polynomials $\mathcal{A}^{\pm}_{U^\tau}(x,y)$ for $U^\tau$ is given by 
\begin{align} \label{formula-AUtaupm}
 \mathcal{A}_{U^\tau}^+(x,y)&=y^2-1-(-1)^\tau xy^{2\tau+2},\\\nonumber  \mathcal{A}^{-}_{U^\tau}(x,y)&=y^2-1+(-1)^{\tau}xy^{2\tau}.  
\end{align}

\begin{remark}
In fact, the dual $A$-polynomial $\mathcal{A}_U(x,y,a)$ and extremal $A$-polynomials $\mathcal{A}_U^\pm (x,y)$ for unknot $U$ has been computed in \cite{GKS} as follow
\begin{align}
    \mathcal{A}_U(x,y,a)&=(y^2-1)-a^{-\frac{1}{2}}x(ay^2-1) \\ \nonumber
    \mathcal{A}_U^+(x,y)&=y^2-1-xy^2\\ \nonumber
    \mathcal{A}_U^-(x,y)&=y^2-1+x
\end{align}
Therefore, by the framing change formulas (\ref{formula-AdualA}) and (\ref{formula-AdualApm}), we can obtain the expressions (\ref{formula-AUtau}) and (\ref{formula-AUtaupm}) directly.  

\end{remark}

Now, we illustrate the computations of the BPS invariants for framed unknot $U^\tau$ with the help of the expression (\ref{formula-AUtau}) for the  dual $A$-polynomial  $\mathcal{A}_{U^\tau}(x,y,a)$. For brevity, we let $X=(-1)^\tau a^{-\frac{1}{2}}x$, and
$Y=1-y^2$, then the dual $A$-polynomial (\ref{formula-AUtau}) is changed to
the functional equation
\begin{align}  \label{mirrorcurvesimple}
Y=X(1-Y)^\tau(1-a(1-Y)).
\end{align}
In order to solve the above equation, we introduce the following
Lagrangian inversion formula \cite{Stan}.
\begin{lemma} \label{lemma-lagrangian}
Let $\phi(\lambda)$ be an invertible formal power series in
indeterminate $\lambda$. Then the functional equation $Y=X\phi(Y)$
has a unique formal power series solution $Y=Y(X)$. Moreover, if $f$
is a formal power series, then
\begin{align} \label{lagrangeinversion}
f(Y(X))=f(0)+\sum_{n\geq
1}\frac{X^n}{n}\left[\frac{df(\lambda)}{d\lambda}\phi(\lambda)^{n}\right]_{\lambda^{n-1}}
\end{align}
\end{lemma}
\begin{remark}
In the following, we will frequently use the binomial coefficient
$\binom{n}{k}$ for all $n\in \mathbb{Z}$. That means that for $n<0$, we
define $\binom{n}{k}=(-1)^k\binom{-n+k-1}{k}$.
\end{remark}

In our case, we take $ \phi(Y)=(1-Y)^\tau(1-a(1-Y)). $ Let
$f(Y)=1-Y$, by formula (\ref{lagrangeinversion}), we obtain
\begin{align*}
y(X)^2=1-Y(X)=1+\sum_{n\geq 1}\frac{X^n}{n}\sum_{i\geq
0}(-1)^{n+i}\binom{n}{i}\binom{n\tau +i}{n-1}a^i
\end{align*}
since $\phi(\lambda)^n$ has the expansion
\begin{align*}
\phi(\lambda)^n&=(1-\lambda)^{n\tau}(1-a(1-\lambda))^n\\\nonumber
&=\sum_{i\geq 0}\binom{n}{i}(-a)^i(1-\lambda)^{n\tau+i} \\\nonumber
&=\sum_{i,j\geq
0}\binom{n}{i}(-1)^{i+j}\binom{n\tau+i}{j}a^i\lambda^j.
\end{align*}
Furthermore, if we let $f(Y(X))=\log(1-Y(X))$, then
\begin{align*}
\left[\frac{df(\lambda)}{d\lambda}\phi(\lambda)^{n}\right]_{\lambda^{n-1}}&=\sum_{i\geq
0}(-1)^{i}\binom{n}{i}\sum_{j=0}^{n-1}(-1)^{j+1}\binom{n\tau+i}{j}a^{i}\\\nonumber
&=\sum_{i\geq
0}(-1)^{i}\binom{n}{i}(-1)^{n}\binom{n\tau+i-1}{n-1}a^{i}
\end{align*}
where we have used the combinatoric identity:
\begin{align*}
\sum_{j=0}^{n-1}(-1)^{j+1}\binom{m}{j}=(-1)^n\binom{m-1}{n-1}.
\end{align*}
Formula (\ref{lagrangeinversion}) gives
\begin{align*}
\log(y(X))=\frac{1}{2}\log(1-Y(X))=\frac{1}{2}\sum_{n\geq 1}\frac{X^n}{n}\sum_{i\geq
0}(-1)^{n+i}\binom{n}{i}\binom{n\tau+i-1}{n-1}a^{i}.
\end{align*}
i.e.
\begin{align*}
\log(y(x))=\frac{1}{2}\sum_{n\geq 1}\frac{x^n}{n}\sum_{i\geq
0}(-1)^{n\tau+n+i}\binom{n}{i}\binom{n\tau+i-1}{n-1}a^{i-\frac{n}{2}}.
\end{align*}

Therefore, 
\begin{align}
x\frac{\partial_x y(x,a)}{y(x,a)}=\frac{1}{2}\sum_{n\geq 1,i\geq 0}(-1)^{n\tau+n+i}\binom{n}{i}\binom{n\tau+i-1}{n-1}x^na^{i-\frac{n}{2}}.    
\end{align}
We introduce
\makeatletter
\let\@@@alph\@alph
\def\@alph#1{\ifcase#1\or \or $'$\or $''$\fi}\makeatother
\begin{subnumcases}
{c_{n,m}(\tau)=}
(-1)^{n\tau+n+\frac{n+m}{2}}\binom{n}{\frac{n+m}{2}}\binom{n\tau+\frac{n+m}{2}-1}{n-1}, &$m+n=even$ and $|m|\leq n$ \label{eq:a1} 
\nonumber \\\nonumber
0, & otherwise.\label{eq:a2}
\end{subnumcases}
\makeatletter\let\@alph\@@@alph\makeatother

By formula (\ref{formula-bri}), we obtain 
\begin{align} \label{formula-bnm}
b_{n,m}(U^\tau)=\frac{1}{n^2}\sum_{d|n,m}\mu(d)c_{\frac{n}{d},\frac{m}{d}}(\tau).    
\end{align}

Similarly, from the expressions (\ref{formula-AUtaupm}) for the extremal $A$-polynomials of $U^\tau$,  we compute the extremal invariants as follows:
\begin{align} \label{formula-bn+}
b^+_{n}(U^\tau)=\frac{1}{n^2}\sum_{d|n}\mu\left(\frac{n}{d}\right)(-1)^{d\tau}\binom{d(\tau+1)-1}{d-1},    
\end{align}
and
\begin{align} \label{formula-bn-}
b^-_{n}(U^\tau)=\frac{1}{n^2}\sum_{d|n}\mu\left(\frac{n}{d}\right)(-1)^{d(\tau+1)}\binom{d\tau-1}{d-1}.    
\end{align}
It is easy to see that $b^\pm_{n}(U^\tau)=b_{n,\pm n}(U^\tau)$, which means that the extremal invariants corresponding to the top and bottom parts of the BPS invariants.

\subsection{Framed twist knots}
\begin{figure}[!htb] 
\begin{align*} 
\raisebox{-15pt}{
\includegraphics[width=120 pt]{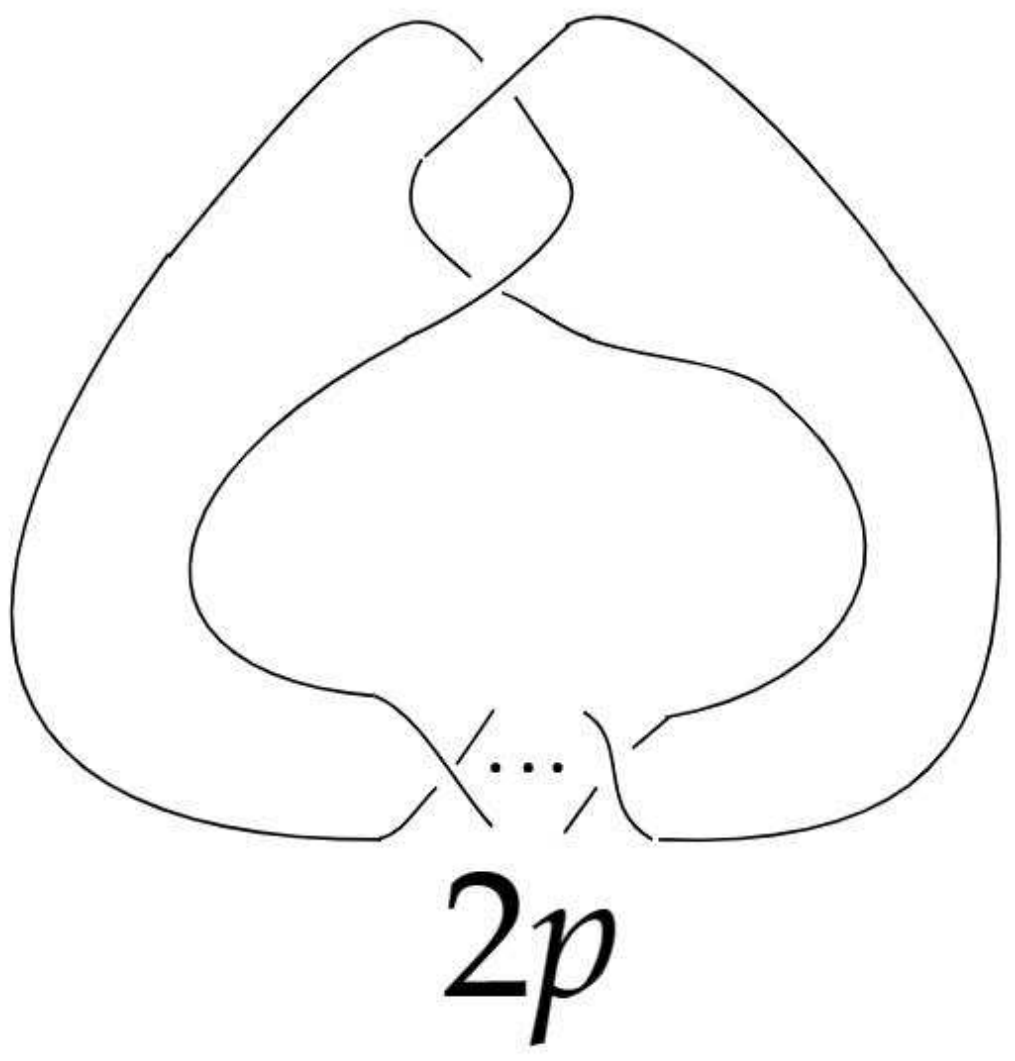}}.
\end{align*}
\caption{Twist knot $\mathcal{K}_p$} \label{figure1} 
\end{figure}

Given an integer $p\in \mathbb{Z}$, the twist knot $\mathcal{K}_p$ is shown in Figure \ref{figure1}. 

Apart from a special case $p=0$ which is the unknot and $p=1$ which is the trefoil knot $3_1$, all other twist knots are hyperbolic. For example, $\mathcal{K}_{-1}=4_1$ and $\mathcal{K}_2=5_2$. The extremal $A$-polynomials for $\mathcal{K}_p$ are given in \cite{GKS}. The formulas for $p\geq 2$ are somewhat different from those for $p\leq -1$. Proposition 1.4 in \cite{GKS} shows that
\begin{align}
\mathcal{A}_{\mathcal{K}_p}^-(x,y)&=x-y^4+y^6, \ \mathcal{A}_{\mathcal{K}_p}^+(x,y)=1-y^2+xy^{4|p|+2}  ,\  p\leq -1; \\\nonumber
\mathcal{A}_{\mathcal{K}_p}^-(x,y)&=1-y^2-xy^4, \ \mathcal{A}_{\mathcal{K}_p}^+(x,y)=1-y^2-xy^{4p+4}  ,\  p\geq 2.
\end{align}

Now we consider the framed twist knot $\mathcal{K}_p^\tau$ with any framing $\tau\in \mathbb{Z}$, by using formula (\ref{formula-AdualApm}), we obtain the explicit formula for the extremal $A$-polynomial of $\mathcal{K}_p^\tau$ as follows.
\begin{align} \label{formula-extremalA1}
\mathcal{A}_{\mathcal{K}_p^\tau}^-(x,y)&=(-1)^\tau xy^{2\tau}-y^4+y^6, \ \ \mathcal{A}_{\mathcal{K}_p^\tau}^+(x,y)=1-y^2+(-1)^\tau xy^{4|p|+2+2\tau} 
\end{align}
for  $p\leq -1$, and 
\begin{align} \label{formula-extremalA2}
\mathcal{A}_{\mathcal{K}_p^\tau}^-(x,y)&=1-y^2-(-1)^\tau x y^{4+2\tau}, \  \ \mathcal{A}_{\mathcal{K}_p^\tau}^+(x,y)=1-y^2-(-1)^\tau xy^{4p+4+2\tau}
\end{align}
for $p\geq 2$.

In the following, we illustrate the computations for the extremal BPS invariant of $\mathcal{K}_p^\tau$ by using the extremal $A$-polynomials (\ref{formula-extremalA1}) and (\ref{formula-extremalA2}).

Let us first consider the case $p\leq -1$. 
We write the equation $\mathcal{A}_{\mathcal{K}_p^\tau}^-(x,y)=0$ as
\begin{align}
    1-y^2=(-1)^\tau xy^{2\tau-4}.
\end{align}
For brevity, let $Y=1-y^2$ and $X=(-1)^\tau x$, then we obtain the equation 
\begin{align}
    Y=X(1-Y)^{\tau-2}. 
\end{align}
We can solve this equation using the Lagrangian inversion formula shown in Lemma \ref{lemma-lagrangian}. 

Let $\phi(Y)=(1-Y)^{\tau-2}$, then
\begin{align}
\phi(\lambda)^n=(1-\lambda)^{n(\tau-2)}=\sum_{j\geq 0}\binom{n(\tau-2)}{j}(-1)^j\lambda^j    
\end{align}
and
\begin{align}
\left[-\frac{1}{1-\lambda}\phi(\lambda)^n\right]_{\lambda^{n-1}}=-\sum_{j=0}^{n-1}\binom{n(\tau-2)}{j}(-1)^j=-\binom{n(3-\tau)-1}{n-1}.    
\end{align}
Then, by formula (\ref{lagrangeinversion}), we obtain
\begin{align}
\log(1-Y(X))=-\sum_{n\geq 1}\frac{X^n}{n}  \binom{n(3-\tau)-1}{n-1}.   
\end{align}
In other words, we have
\begin{align}
\log (y(x))=-\frac{1}{2} \sum_{n\geq 1}\frac{(-1)^{n\tau} x^n}{n}  \binom{n(3-\tau)-1}{n-1}.   
\end{align}
Then
\begin{align}
 x\frac{\partial_x y}{y}=x\frac{d\log (y(x))}{dx}=-\frac{1}{2}\sum_{n\geq 1}(-1)^{n\tau}x^n  \binom{n(3-\tau)-1}{n-1}.  
\end{align}
Finally, using the formula (\ref{formula-bri}), we obtain 
\begin{align} \label{formula-bn-p1}
 b_{r}^-(\mathcal{K}^\tau_p)=-\frac{1}{r^2}\sum_{d|r}\mu\left(\frac{r}{d}\right)(-1)^{d\tau}\binom{d(3-\tau)-1}{d-1}. 
\end{align}

Now we consider the equation 
\begin{align}
  \mathcal{A}_{\mathcal{K}_p^\tau}^+(x,y)=1-y^2+(-1)^\tau xy^{4|p|+2+2\tau}=0.   
\end{align}
We rewrite it as 
\begin{align}
Y=X(1-Y)^{2|p|+\tau+1}  
\end{align}
where $X=(-1)^{\tau-1}x$ and $Y=1-y^2$. We solve it using the Lagrangian inversion formula (\ref{lagrangeinversion}). Taking $\phi(Y)=(1-Y)^{2|p|+\tau+1}$, we have

\begin{align}
\log(1-Y(X))&=\sum_{n\geq 1}\frac{X^n}{n}\left[-\frac{1}{1-\lambda}\phi(\lambda)^n\right]_{\lambda^{n-1}}\\\nonumber
&=\sum_{n\geq 1}\frac{(-1)^{\tau(n-1)}x^n}{n}(-1)^{n}\binom{n(2|p|+1+\tau)-1}{n-1}.   
\end{align}
Then we have
\begin{align}
x\frac{\partial_xy}{y}=x\frac{d\log(y(x))}{dx}=\frac{1}{2}\sum_{n\geq 1}(-1)^{n\tau}x^n \binom{n(2|p|+1+\tau)-1}{n-1}.   
\end{align}
Therefore, using (\ref{formula-bri}),  we obtain
\begin{align}   \label{formula-bn+p1}
b_{r}^+(\mathcal{K}^\tau_p)=\frac{1}{r^2}\sum_{d|r}\mu\left(\frac{r}{d}\right)(-1)^{d\tau}\binom{d(2|p|+1+\tau)-1}{d-1}.     
\end{align}
For the case $p\geq 2$, with the similar computations shown above, it is easy to obtain 
\begin{align}  \label{formula-bn-p2}
 b_r^-(\mathcal{K}^\tau_p)=\frac{1}{r^2}\sum_{d|r}\mu\left(\frac{r}{d}\right)(-1)^{d(\tau+1)}\binom{d(\tau+2)-1}{d-1}  
\end{align}
and
\begin{align} \label{formula-bn+p2}
 b_{r}^+(\mathcal{K}^\tau_p)=\frac{1}{r^2}\sum_{d|r}\mu\left(\frac{r}{d}\right)(-1)^{d(\tau+1)}\binom{d(\tau+2+2p)-1}{d-1}.   
\end{align}

\subsection{Integrality} 
Based on the integrality conjecture for the Ooguri-Vafa invariants and BPS invariants, the above formulas (\ref{formula-bnm}), (\ref{formula-bn+}), (\ref{formula-bn-}), (\ref{formula-bn-p1}), (\ref{formula-bn+p1}), (\ref{formula-bn-p2}) and (\ref{formula-bn+p2}) 
will
give integers. However,  the integrality of these formulas is not so obvious from their expressions. A straightforward proof is required. 

In fact, in the previous work of the second author and W. Luo \cite{LZhu19}, we have proved that formula (\ref{formula-bnm}) really gives integers. See Theorem 1.1 in \cite{LZhu19}.  Moreover, in \cite{Zhu19}, the second author proved the following.
\begin{theorem}(Theorem 5.3 in \cite{Zhu19})\label{theorem-integrality}
For any $r\geq 1$ and $t \in \mathbb{Z}$, then we have
\begin{align} \label{formula-integrality}
 \frac{1}{r^2}\sum_{d|r}\mu\left(\frac{r}{d}\right)(-1)^{d (t+1)}\binom{dt-1}{d-1}\in \mathbb{Z}   
\end{align}
\end{theorem}
Then the integrality of $b_r^{\pm}(U^\tau)$, $b_{r}^\pm(\mathcal{K}_p^\tau)$ is an easy consequence of Theorem \ref{theorem-integrality}.  For example, by letting $t=\tau+2+2p$, then the formula (\ref{formula-integrality}) gives 
\begin{align}
     b_{r}^+(\mathcal{K}^\tau_p)=\frac{1}{r^2}\sum_{d|r}\mu\left(\frac{r}{d}\right)(-1)^{d(\tau+1)}\binom{d(\tau+2+2p)-1}{d-1}\in \mathbb{Z}. 
\end{align}

\section{BPS invariants of framed links} \label{section-link}
Given a framed link $\mathcal{L}$ with the $L$ components $\mathcal{K}_1,..,\mathcal{K}_L$. Suppose that the framing of $\mathcal{K}_\alpha$ is given by $\tau_\alpha$, we let $\vec{\tau}=(\tau_1,...,\tau_L)$ and denote the framed link $\mathcal{L}$ by $\mathcal{L}^{\vec{\tau}}$. 
The framed colored HOMFLYPT invariant (with symmetric representations) for $\mathcal{L}^{\vec{\tau}}$ is defined by 
\begin{align}
\mathcal{H}_{r_1,...,r_L}(\mathcal{L}^{\vec{\tau}};q,a)=(-1)^{\sum_{\alpha=1}^Lr_\alpha|\tau_\alpha|}q^{\sum_{\alpha=1}^L\frac{r_\alpha(r_\alpha-1)\tau_\alpha}{2}}H_{r_1,...,r_L}(\mathcal{L};q,a),  
\end{align}
where $H_{r_1,...,r_L}(\mathcal{L};q,a)$ denotes the framing-independent colored HOMFLYPT invariant with each component $\mathcal{K}_\alpha$ colored by symmetric representation labelled by  $r_\alpha\in \mathbb{N}$.  

We consider the following generating function for  $\mathcal{H}_{r_1,...,r_L}(\mathcal{L}^{\vec{\tau}};q,a)$, 
\begin{align} \label{formula-Zlink}
\mathcal{Z}(\mathcal{L}^{\vec{\tau}};q,a;x_{1},x_{2},...,x_L)=\sum_{r_1,...,r_L\geq
0}\mathcal{H}_{r_1,...,r_L}(\mathcal{L}^{\vec{\tau}};q,a)x_1^{r_1}\cdots x_L^{r_{L}}.
\end{align}
Parallel to formula (\ref{formula-reducedCSknot}), the generating function can be rewritten in the following form 
\begin{align}  
    \mathcal{Z}(\mathcal{L}^{\vec{\tau}};q,a;x_{1},...,x_L)=\exp\left(\sum_{d\geq 1}\sum_{\substack{r_1,r_2,...,r_L\geq 0\\ (r_1,..,r_L)\neq (0,..,0)}}\frac{1}{d}\mathfrak{f}_{r_1,...,r_L}(\mathcal{L}^{\vec{\tau}};a^d,q^d)x_1^{dr_1}\cdots x_L^{dr_L})\right). 
\end{align}

By M\"obius inversion formula, we obtain
\begin{align}
\mathfrak{f}_{r_1,...,r_L}(\mathcal{L}^{\vec{\tau}};q,a) =[x_1^{r_1}\cdots
x_L^{r_L}]\sum_{d\geq 1}\frac{\mu(d)}{d} \Psi_d(\log(\mathcal{Z}(\mathcal{L}^{\vec{\tau}};q,a;x_{1},...,x_L))),
\end{align}
where the notation $[x_1^{r_1}\cdots x_L^{r_L}]f(x_1,...,x_L,q,a)$
denotes the coefficients of $x_1^{r_1}\cdots x_L^{r_L}$ in the
function $f(x_1,...,x_L,q,a)$ and the operator $\Psi_d$ is defined as $\Psi_d(f(x_1,...,x_L,q,a)=f(x_1^d,..,x_L^d,q^d,a^d)$.

We can generalize the original Ooguri-Vafa integrality conjecture \cite{OV} to the case of framed link $\mathcal{L}^{\vec{\tau}}$ as follows.  

\begin{conjecture}
Suppose $r_1,r_2,...,r_L\geq 0$ and $(r_1,..,r_L)\neq (0,..,0)$, let $k=\#\{r_i|r_i\neq 0, i=1,...,L\}$, and hence $1\leq k\leq L$. There exists integers $N_{(r_1,...,r_L),i,j}(\mathcal{L}^{\vec{\tau}})$ which vanish for large $|i|,|j|$, such that  
\begin{align}
    \mathfrak{f}_{r_1,...,r_L}(\mathcal{L}^{\vec{\tau}};q,a)=(q^{\frac{1}{2}}-q^{-\frac{1}{2}})^{k-2}\sum_{i,j\in \mathbb{Z}}N_{(r_1,...,r_L),i,j}(\mathcal{L}^{\vec{\tau}})a^{\frac{i}{2}}q^{\frac{j}{2}}.
\end{align}
\end{conjecture}

\begin{remark}
Using the strong integerality property of the colored HOMFLYPT invariants proved in \cite{Zhu23}, the above conjecture can be refined as follows. 

Under the above setting, we have
\begin{align}
\mathfrak{f}_{r_1,...,r_L}(\mathcal{L}^{\vec{\tau}};q,a)=(q^{\frac{1}{2}}-q^{-\frac{1}{2}})^{k-2}
p_{r_1,...,r_L}(\mathcal{L}^{\vec{\tau}};q,a),
\end{align}
and where
\begin{align}
p_{r_1,...r_L}(\mathcal{L}^{\vec{\tau}};q,a)=\sum_{i,j\in \mathbb{Z}}N_{(r_1,...,r_L),i,j}(\mathcal{L}^{\vec{\tau}})a^{\frac{i}{2}}q^{\frac{j}{2}}\in
a^{\frac{1}{2}\epsilon_1}q^{\frac{1}{2}\epsilon_2} \mathbb{Z}[a^{\pm 1},q^{\pm 1}]
\end{align}
where $\epsilon_1,\epsilon_2\in \{0,1\}$ is determined by $$ |r_1+\cdots +r_L|\equiv
\epsilon_1 \mod 2 $$ and $$ |r_1+\cdots +r_L|+k\equiv
\epsilon_2 \mod 2. $$ 
\end{remark}

We define the BPS invariant for $\mathcal{L}^{\vec{\tau}}$ by
\begin{align}
b_{(r_1,...,r_L),i}(\mathcal{L}^{\vec{\tau}})=\sum_{j\in \mathbb{Z}} N_{(r_1,...,r_L),i,j}(\mathcal{L}^{\vec{\tau}}).    
\end{align}
The generating function of these BPS invariants is given by
\begin{align}
    p_{r_1,...r_L}(\mathcal{L}^{\vec{\tau}};1,a)=a^{\frac{\epsilon}{2}}\sum_{i\in \mathbb{Z}}b_{(r_1,...,r_L),i}(\mathcal{L}^{\vec{\tau}})a^{i}.
\end{align}
Then it is conjectured that
\begin{align}
    p_{r_1,...r_L}(\mathcal{L}^{\vec{\tau}};1,a)\in a^{\frac{1}{2}\epsilon_1} \mathbb{Z}[a^{\pm 1}] 
\end{align}
which implies the integrality of the BPS invariants $b_{(r_1,...,r_L),i}(\mathcal{L}^{\vec{\tau}})$.

\subsection{General formulas}
Let $\vec{r}=(r_1,...,r_L)\in (\mathbb{Z}_{\geq 0})^L$, then the set $S=\{\vec{r}|\vec{r}\in (\mathbb{Z}_{\geq 0})^L, \vec{r}\neq (0,...,0) \}$ has a total order in lexicographical order. Suppose that it is given by $\vec{r}^{(1)}<\vec{r}^{(2)}<\vec{r}^{(3)}<\cdots <\vec{r}^{(n)} <\cdots$. 
Given a sequence of nonnegative integers $k_i\in \mathbb{Z}_{\geq 0}$ with $i=1,2,...$, then
$\mathfrak{U}=((\vec{r}^{(1)})^{k_1}(\vec{r}^{(2)})^{k_2}\cdots )$ forms a partition of the vector $\vec{r}=\sum_{i\geq 1}k_i\vec{r}^{(i)}$. We introduce the notation $|\mathfrak{U}|=\vec{r}$, $l(\mathfrak{U})=\sum_{i}k_i$, $|Aut(\mathfrak{U})|=\prod_{i\geq 1}(k_i)!$ and
$\mathcal{H}_{\mathfrak{U}}(\mathcal{L}^{\vec{\tau}};q,a)=\prod_{i\geq 1}(\mathcal{H}_{\vec{r}^{(i)}}(\mathcal{L}^{\vec{\tau}};q,a))^{k_i}$.
Then the formula (\ref{formula-Zlink}) can be written as
\begin{align}
\mathcal{Z}(\mathcal{L}^{\vec{\tau}};q,a,\vec{x})=1+\sum_{\vec{r}\in S}\mathcal{H}_{\vec{r}}(\mathcal{L}^{\vec{\tau}};q,a)\vec{x}^{\vec{r}}.  
\end{align}
We make the following expansion
\begin{align}
    \log(\mathcal{Z}(\mathcal{L}^{\vec{\tau}};q,a,\vec{x}))=\sum_{\vec{r}\in S}\mathcal{F}_{\vec{r}}(\mathcal{L}^{\vec{\tau}};q,a)\vec{x}^{\vec{r}}, 
\end{align}
where
\begin{align} \label{formula-mathcalF}
  \mathcal{F}_{\vec{r}}(\mathcal{L}^{\vec{\tau}};q,a)=\sum_{|\mathfrak{U}|=\vec{r}}(-1)^{l(\mathfrak{U})}\frac{(l(\mathfrak{U})-1)!}{|Aut(\mathfrak{U})|}\mathcal{H}_{\mathfrak{U}}(\mathcal{L}^{\vec{\tau}};q,a).  
\end{align}
Therefore, we obtain
\begin{align}
p_{\vec{r}}(\mathcal{L}^{\vec{\tau}};q,a)=(q^{\frac{1}{2}}-q^{-\frac{1}{2}})^{2-k}\sum_{d|\vec{r}}\frac{\mu(d)}{d}\mathcal{F}_{\frac{\vec{r}}{d}}(\mathcal{L}^{\vec{\tau}};q^d,a^d),
\end{align}
and
\begin{align}
p_{\vec{r}}(\mathcal{L}^{\vec{\tau}};1,a)=\lim_{q\rightarrow 1}(q^{\frac{1}{2}}-q^{-\frac{1}{2}})^{2-k}\sum_{d|\vec{r}}\frac{\mu(d)}{d}\mathcal{F}_{\frac{\vec{r}}{d}}(\mathcal{L}^{\vec{\tau}};q^d,a^d).    
\end{align}

\subsection{Framed Whitehead link $W^{(\tau_1,\tau_2)}$}

For an integer $n$, we define the symbols by
\begin{align*}
[n] &= \frac{q^{\frac{n}{2}}-q^{-\frac{n}{2}}}{q^{\frac{1}{2}}-q^{-\frac{1}{2}}}, \qquad \{ n \} =  q^{\frac{n}{2}} -
q^{-\frac{n}{2}}, \qquad \{ n ; a \} =  a^{\frac{1}{2}}q^{\frac{n}{2}}-a^{-\frac{1}{2}}q^{-\frac{n}{2}}.
\end{align*}
For integers $n >0 , i \geq 0$, we introduce the products of  $i$
terms of these symbols by
\begin{align*}
[n]_{i} &= [n] [n-1] \cdots [n-i+1], \\
\{ n \}_{i} &=  \{ n \} \{ n-1 \} \cdots \{ n-i+1 \},  \\
\{ n ; a \}_{i} &=  \{ n; a \} \{ n-1; a \} \cdots \{ n-i+1; a \},  \\
\{ -n ; a \}_{i} &=  \{ -n; a \} \{ -n+1; a \} \cdots \{ -n+i-1; a
\},
\end{align*}
which are defined to be $1$ if $i = 0$. Furthermore, we let
\begin{align*}
[n]! =[n]_{n}, \quad \{ n \} ! =  \{ n \}_{n}, \quad \left[
\begin{array}{@{\,}c@{\,}}n \\ i \end{array} \right] =
\frac{[n]!}{[i]! [n-i]!}.
\end{align*}

\begin{figure}[!htb] 
\begin{align*}
\raisebox{-50pt}{
\includegraphics[width=140 pt]{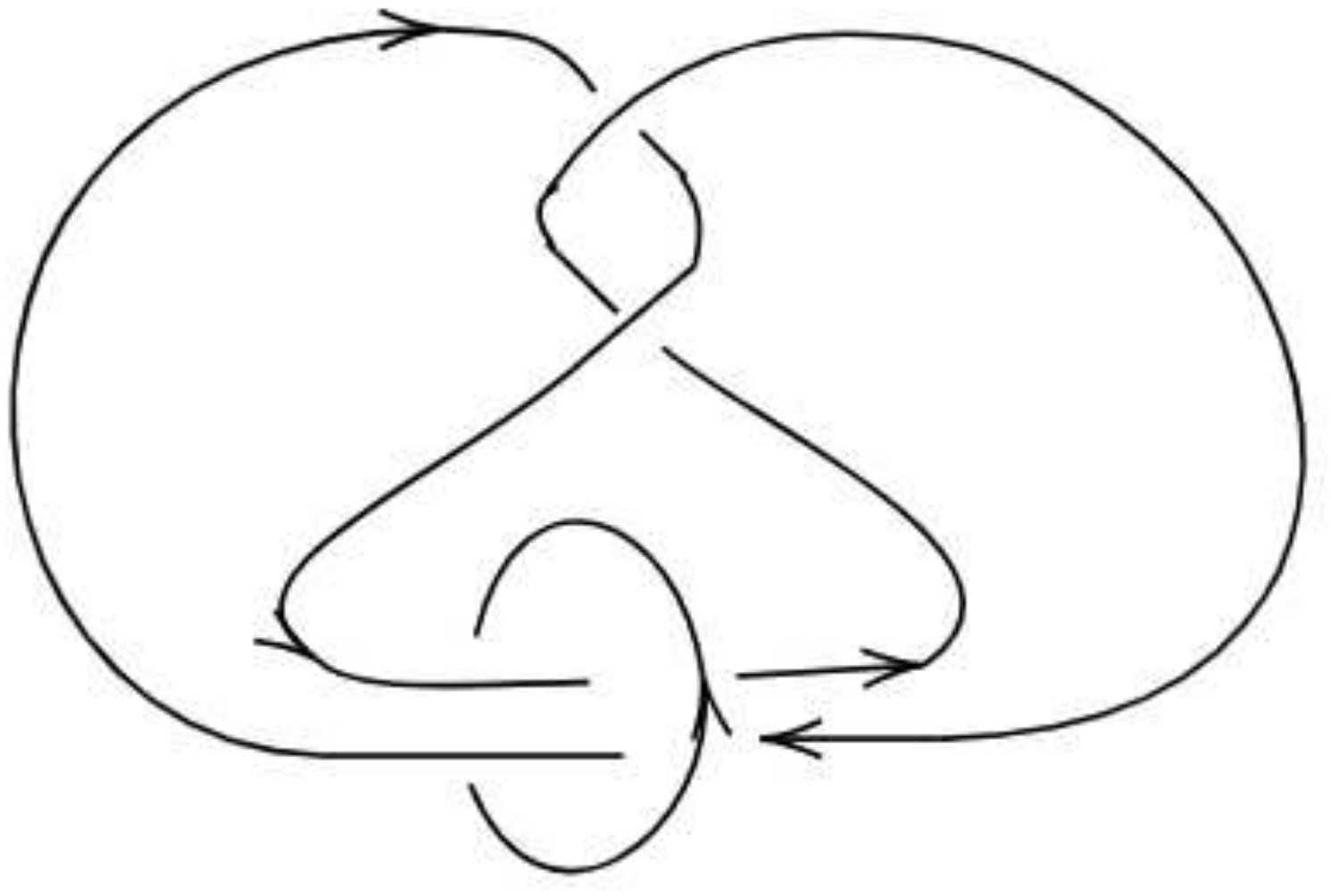}}
\end{align*}
\caption{Whitehead link $W$} \label{fig:WH}
\end{figure}

The Whitehead link $W$ is a 2-component link shown in Figure \ref{fig:WH}. 

The colored HOMFLYPT invariant of $W$ was predicted in \cite{GNSSS15} and computed in \cite{CLZ} using the HOMFLYPT skein theory.  
\begin{align} \label{formula-WL}
H_{r_1,r_2}(W;q,a)
=&\sum_{i=0}^{\min{\{r_1,r_2\}}}(-1)^i \frac{ \{ r_1+i-1;a
\}_{r_1-i} }{ \{ r_1-i\}! }\frac{ \{ r_2+i-1;a \}_{r_2-i} }{ \{
r_2-i\}! }\\\nonumber &\cdot
\frac{a^{\frac{i}{2}}q^{\frac{i(i-1)}{4}}}{\{i\}!}\{2i-1;a\}_{2i}\{i-2;a\}_i\{i\}!.
\end{align}
then
\begin{align}
\mathcal{H}_{r_1,r_2}(W^{(\tau_1,\tau_2)};q,a)=(-1)^{r_1\tau_1+r_2\tau_2}q^{\frac{r_1(r_1-1)\tau_1+r_2(r_2-1)\tau_2}{2}}H_{r_1,r_2}(W;q,a).     
\end{align}
By formula (\ref{formula-mathcalF}), some lower-degree terms are given as follows. 
\begin{align*}
\mathcal{F}_{1,0}(q,a)=\mathcal{H}_{1,0}(q,a),  \ \ \mathcal{F}_{0,1}(q,a)=\mathcal{H}_{0,1}(q,a),   
\end{align*}
\begin{align*}
\mathcal{F}_{2,0}(q,a)=\mathcal{H}_{2,0}(q,a)-\frac{1}{2}\mathcal{H}_{1,0}(q,a)^2,  \ \ \mathcal{F}_{0,2}(q,a)=\mathcal{H}_{0,2}(q,a)-\frac{1}{2}\mathcal{H}_{0,1}(q,a)^2,   
\end{align*}
\begin{align*}
\mathcal{F}_{1,1}(q,a)=\mathcal{H}_{1,1}(q,a)-\mathcal{H}_{1,0}(q,a)\mathcal{H}_{0,1}(q,a),   
\end{align*}
\begin{align*}
&\mathcal{F}_{2,1}=\mathcal{H}_{2,1}-\mathcal{H}_{0,1}\mathcal{H}_{2,0}-\mathcal{H}_{1,0}\mathcal{H}_{1,1}+\mathcal{H}_{0,1}\mathcal{H}_{1,0}^2,   
\end{align*}
\begin{align*}
&\mathcal{F}_{1,2}=\mathcal{H}_{1,2}-\mathcal{H}_{1,0}\mathcal{H}_{0,2}-\mathcal{H}_{0,1}\mathcal{H}_{1,1}+\mathcal{H}_{1,0}\mathcal{H}_{0,1}^2,   
\end{align*}
\begin{align*}
&\mathcal{F}_{2,2}=\mathcal{H}_{2,2}-\mathcal{H}_{1,0}\mathcal{H}_{1,2}-\frac{3}{2}\mathcal{H}_{1,0}^2\mathcal{H}_{0,1}^2+\mathcal{H}_{1,0}^2\mathcal{H}_{0,2}\\\nonumber
&+2\mathcal{H}_{1,0}\mathcal{H}_{0,1}\mathcal{H}_{1,1}-\mathcal{H}_{0,1}\mathcal{H}_{2,1}+\mathcal{H}_{0,1}^2\mathcal{H}_{2,0}-\mathcal{H}_{2,0}\mathcal{H}_{0,2}-\frac{1}{2}\mathcal{H}_{1,1}^2.
\end{align*}

By straightforward computations, we obtain 
\begin{align*}
p_{1,1}(W^{(\tau_1,\tau_2)},1,a)=-(-1)^{\tau_1+\tau_2}a^2+3(-1)^{\tau_1+\tau_2}a-3(-1)^{\tau_1+\tau_2}+(-1)^{\tau_1+\tau_2}a^{-1},   
\end{align*}
\begin{align*}
 p_{2,0}(W^{(\tau_1,\tau_2)},1,a)=\frac{1-(-1)^{\tau_1}+2\tau_1}{4}a-\tau_1+\frac{-1+(-1)^{\tau_1}+2\tau_1}{4}a^{-1},   
\end{align*}
\begin{align*}
  p_{2,1}(W^{(\tau_1,\tau_2)},1,a)&=((-1)^{\tau_2}(1+\tau_1))a^{\frac{5}{2}}-(4(-1)^{\tau_2}\tau_1+2(-1)^{\tau_2})a^{\frac{3}{2}}+6\tau_1(-1)^{\tau_2}a^{\frac{1}{2}}\\\nonumber
  &+(2(-1)^{\tau_2}-4(-1)^{\tau_2}\tau_1)a^{-\frac{1}{2}}+(-1)^{\tau_2}(\tau_1-1)a^{-\frac{3}{2}}, 
\end{align*}
\begin{align*}
p_{2,2}(W^{(\tau_1,\tau_2)},1,a)&=\frac{1+(-1)^{\tau_1+\tau_2}}{2}a^4-(\tau_1\tau_2+\tau_1+\tau_2+4)a^3\\\nonumber
&+(5\tau_1\tau_2+3\tau_1+3\tau_2+\frac{17-3(-1)^{\tau_1+\tau_2}}{2})a^2\\\nonumber
&-(2\tau_1-2\tau_2-10\tau_1\tau_2-8)a-2\tau_1-2\tau_2+10\tau_1\tau_2+\frac{11+3(-1)^{\tau_1+\tau_2}}{2}\\\nonumber
&+(3\tau_1+3\tau_2-5\tau_1\tau_2-4)a^{-1}+(\tau_1\tau_2-\tau_1-\tau_2+\frac{3-(-1)^{\tau_1+\tau_2}}{2})a^{-2}.
\end{align*}

These formulas confirm the integrality of the BPS invariants for framed Whitehead links $W^{(\tau_1,\tau_2)}$. 

Furthermore, we list some Ooguri-Vafa invariants in the following tables, which confirm the Conjecture \ref{conjecture-framedlink}.  

\begin{table}[H]
\makebox[\textwidth]{
\centering
\renewcommand{\arraystretch}{1.25}
$\begin{array}{|c|cccccccc|}
\hline
i \backslash j & 4 & 3 & 2 & 1 & 0 & -1 & -2 & -3 \\ \hline
4 & 0 & 1 & 0 & 0 & 0 & 0 & 0 & 0 \\
3 & -1 & -2 & -1 & -1 & 0 & 1 & 0 & 0 \\
2 & 2 & 2 & 3 & 2 & -1 & -1 & 0 & 0 \\
1 & -1 & -2 & -3 & -1 & 1 & 0 & -1 & -1 \\
0 & 0 & 1 & 1 & 0 & 0 & 0 & 3 & 2 \\
-1 & 0 & 0 & 0 & 0 & 1 & -1 & -3 & -1 \\
-2 & 0 & 0 & 0 & 0 & -1 & 1 & 1 & 0 \\ \hline
\end{array}$
}
\newline \\ 
\caption{Ooguri-Vafa invariants $N_{(2,2),i,j}(W^{(0,0)})$}
\label{}
\end{table}

\begin{table}[H]
\makebox[\textwidth]{
\centering
\renewcommand{\arraystretch}{1.25}
$\begin{array}{|c|cccccccc|}
\hline
i \backslash j & 5 & 4 & 3 & 2 & 1 & 0 & -1 & -2 \\ \hline
4 & 0 & 1 & 0 & 0 & 0 & -1 & 0 & 0 \\
3 & -1 & -2 & -1 & -2 & -1 & 1 & 1 & 0 \\
2 & 2 & 2 & 4 & 5 & 2 & -1 & -1 & 0 \\
1 & -1 & -2 & -4 & -4 & -2 & 2 & 1 & 0 \\
0 & 0 & 1 & 1 & 2 & 0 & -2 & -1 & 1 \\
-1 & 0 & 0 & 0 & -1 & 1 & 1 & -1 & -1 \\
-2 & 0 & 0 & 0 & 0 & 0 & 0 & 1 & 0 \\ \hline
\end{array}$
}
\newline \\ 
\caption{Ooguri-Vafa invariants $N_{(2,2),i,j}(W^{(0,1)})$}
\label{}
\end{table}

\begin{table}[H]
\makebox[\textwidth]{
\centering
\renewcommand{\arraystretch}{1.25}
$\begin{array}{|c|cccccccc|}
\hline
i \backslash j & 5 & 4 & 3 & 2 & 1 & 0 & -1 & -2 \\ \hline
4 & 0 & 1 & 0 & 0 & 0 & -1 & 0 & 0 \\
3 & -1 & -2 & -1 & -2 & -1 & 1 & 1 & 0 \\
2 & 2 & 2 & 4 & 5 & 2 & -1 & -1 & 0 \\
1 & -1 & -2 & -4 & -4 & -2 & 2 & 1 & 0 \\
0 & 0 & 1 & 1 & 2 & 0 & -2 & -1 & 1 \\
-1 & 0 & 0 & 0 & -1 & 1 & 1 & -1 & -1 \\
-2 & 0 & 0 & 0 & 0 & 0 & 0 & 1 & 0 \\ \hline
\end{array}$
}
\newline \\ 
\caption{Ooguri-Vafa invariants $N_{(2,2),i,j}(W^{(1,0)})$}
\label{}
\end{table}

\begin{table}[H]
\makebox[\textwidth]{
\centering
\renewcommand{\arraystretch}{1.25}
$\begin{array}{|c|ccccccccc|}
\hline
i \backslash j & 6 & 5 & 4 & 3 & 2 & 1 & 0 & -1 & -2 \\ \hline
4 & 0 & 1 & 0 & 0 & 0 & 0 & 0 & 0 & 0 \\
3 & -1 & -2 & -1 & -2 & -3 & 0 & 1 & 1 & 0 \\
2 & 2 & 2 & 4 & 6 & 6 & 2 & -1 & -2 & -1 \\
1 & -1 & -2 & -4 & -7 & -7 & -4 & 1 & 1 & 1 \\
0 & 0 & 1 & 1 & 4 & 5 & 2 & 0 & 0 & 0 \\
-1 & 0 & 0 & 0 & -1 & -1 & 0 & -1 & 0 & 0 \\ \hline
\end{array}$
}
\newline \\ 
\caption{Ooguri-Vafa invariants $N_{(2,2),i,j}(W^{(1,1)})$}
\label{}
\end{table}

\begin{table}[H]
\makebox[\textwidth]{
\centering
\renewcommand{\arraystretch}{1.25}
$\begin{array}{|c|ccccccccccc|}
\hline
a^i \backslash q^j & \frac{11}{2} & \frac{9}{2} & \frac{7}{2} & \frac{5}{2} & \frac{3}{2} &
\frac{1}{2} & -\frac{1}{2} & -\frac{3}{2} & -\frac{5}{2} & -\frac{7}{2} &
-\frac{9}{2} \\ \hline
\frac{9}{2} & 0 & 1 & 1 & 0 & 0 & -1 & 0 & 0 & 0 & 0 & 0 \\
\frac{7}{2} & -1 & -3 & -3 & -1 & 0 & 1 & 1 & 0 & 0 & 0 & 0 \\
\frac{5}{2} & 2 & 4 & 4 & 3 & 0 & -2 & -2 & 1 & 0 & 0 & 0 \\
\frac{3}{2} & -1 & -3 & -3 & -2 & 1 & 4 & 0 & -1 & 0 & 1 & 0 \\
\frac{1}{2} & 0 & 1 & 1 & 0 & -2 & -2 & 1 & 1 & -1 & -3 & -1 \\
-\frac{1}{2} & 0 & 0 & 0 & 0 & 1 & 0 & -1 & -2 & 4 & 4 & 2 \\
-\frac{3}{2} & 0 & 0 & 0 & 0 & 0 & 0 & 2 & 0 & -4 & -3 & -1 \\
-\frac{5}{2} & 0 & 0 & 0 & 0 & 0 & 0 & -1 & 1 & 1 & 1 & 0 \\ \hline
\end{array}$
}
\newline \\ 
\caption{Ooguri-Vafa invariants $N_{(2,3),i,j}(W^{(0,0)})$}
\label{}
\end{table}

\begin{table}[H]
\makebox[\textwidth]{
\centering
\renewcommand{\arraystretch}{1.25}
$\begin{array}{|c|ccccccccccc|}
\hline
i \backslash j & \frac{17}{2} & \frac{15}{2} & \frac{13}{2} & \frac{11}{2} & \frac{9}{2} &
\frac{7}{2} & \frac{5}{2} & \frac{3}{2} & \frac{1}{2} & -\frac{1}{2} &
-\frac{3}{2} \\ \hline
\frac{9}{2} & 0 & -1 & -1 & -1 & -1 & 0 & 0 & 1 & 1 & 0 & 0 \\
\frac{7}{2} & 1 & 3 & 4 & 5 & 6 & 4 & 2 & -2 & -4 & -2 & 0 \\
\frac{5}{2} & -2 & -4 & -7 & -11 & -13 & -12 & -7 & 2 & 7 & 4 & 1 \\
\frac{3}{2} & 1 & 3 & 6 & 11 & 14 & 15 & 9 & -2 & -8 & -4 & -1 \\
\frac{1}{2} & 0 & -1 & -2 & -5 & -8 & -9 & -5 & 2 & 6 & 3 & -1 \\
-\frac{1}{2} & 0 & 0 & 0 & 1 & 2 & 2 & 2 & -2 & -3 & 0 & 1 \\
-\frac{3}{2} & 0 & 0 & 0 & 0 & 0 & 0 & -1 & 1 & 1 & -1 & 0 \\ \hline
\end{array}$
}
\newline \\ 
\caption{Ooguri-Vafa invariants $N_{(2,3),i,j}(W^{(0,1)})$}
\label{}
\end{table}

\begin{table}[H]
\makebox[\textwidth]{
\centering
\renewcommand{\arraystretch}{1.25}
$\begin{array}{|c|ccccccccccc|}
\hline
i \backslash j & \frac{13}{2} & \frac{11}{2} & \frac{9}{2} & \frac{7}{2} & \frac{5}{2} &
\frac{3}{2} & \frac{1}{2} & -\frac{1}{2} & -\frac{3}{2} & -\frac{5}{2} &
-\frac{7}{2} \\ \hline
\frac{9}{2} & 0 & 1 & 1 & 0 & 0 & 0 & -1 & 0 & 0 & 0 & 0 \\
\frac{7}{2} & -1 & -3 & -3 & -1 & -2 & 0 & 2 & 1 & 0 & 0 & 0 \\
\frac{5}{2} & 2 & 4 & 4 & 4 & 4 & -1 & -3 & -2 & 1 & 0 & 0 \\
\frac{3}{2} & -1 & -3 & -3 & -4 & -2 & 2 & 5 & 0 & -1 & 0 & 0 \\
\frac{1}{2} & 0 & 1 & 1 & 1 & 0 & -3 & -4 & 2 & 1 & -1 & -1 \\
-\frac{1}{2} & 0 & 0 & 0 & 0 & 0 & 3 & 0 & -2 & -1 & 3 & 2 \\
-\frac{3}{2} & 0 & 0 & 0 & 0 & 0 & -1 & 1 & 2 & -1 & -3 & -1 \\
-\frac{5}{2} & 0 & 0 & 0 & 0 & 0 & 0 & 0 & -1 & 1 & 1 & 0 \\  \hline
\end{array}$
}
\newline \\ 
\caption{Ooguri-Vafa invariants $N_{(2,3),i,j}(W^{(1,0)})$}
\label{}
\end{table}

\begin{table}[H]
\makebox[\textwidth]{
\centering
\renewcommand{\arraystretch}{1.25}
$\begin{array}{|c|cccccccccccc|}
\hline
i \backslash j & \frac{19}{2} & \frac{17}{2} & \frac{15}{2} & \frac{13}{2} & \frac{11}{2} &
\frac{9}{2} & \frac{7}{2} & \frac{5}{2} & \frac{3}{2} & \frac{1}{2} & -\frac{1}{2}
& -\frac{3}{2} \\ \hline
\frac{9}{2} & 0 & -1 & -1 & -1 & -1 & 0 & 0 & 0 & 1 & 1 & 0 & 0 \\
\frac{7}{2} & 1 & 3 & 4 & 5 & 6 & 5 & 5 & 2 & -3 & -4 & -2 & 0 \\
\frac{5}{2} & -2 & -4 & -7 & -11 & -14 & -16 & -16 & -9 & 2 & 7 & 5 & 1 \\
\frac{3}{2} & 1 & 3 & 6 & 11 & 16 & 21 & 22 & 15 & 1 & -7 & -5 & -2 \\
\frac{1}{2} & 0 & -1 & -2 & -5 & -9 & -13 & -16 & -12 & -2 & 4 & 3 & 1 \\
-\frac{1}{2} & 0 & 0 & 0 & 1 & 2 & 3 & 6 & 5 & 1 & -1 & -1 & 0 \\
-\frac{3}{2} & 0 & 0 & 0 & 0 & 0 & 0 & -1 & -1 & 0 & 0 & 0 & 0 \\ \hline
\end{array}$
}
\newline \\ 
\caption{Ooguri-Vafa invariants $N_{(2,3),i,j}(W^{(1,1)})$}
\label{}
\end{table}

\subsection{Framed Borromean rings $\mathbf{B}^{(\tau_1,\tau_2,\tau_3)}$}

\begin{figure}
\begin{align*}
\raisebox{-50pt}{
\includegraphics[width=140 pt]{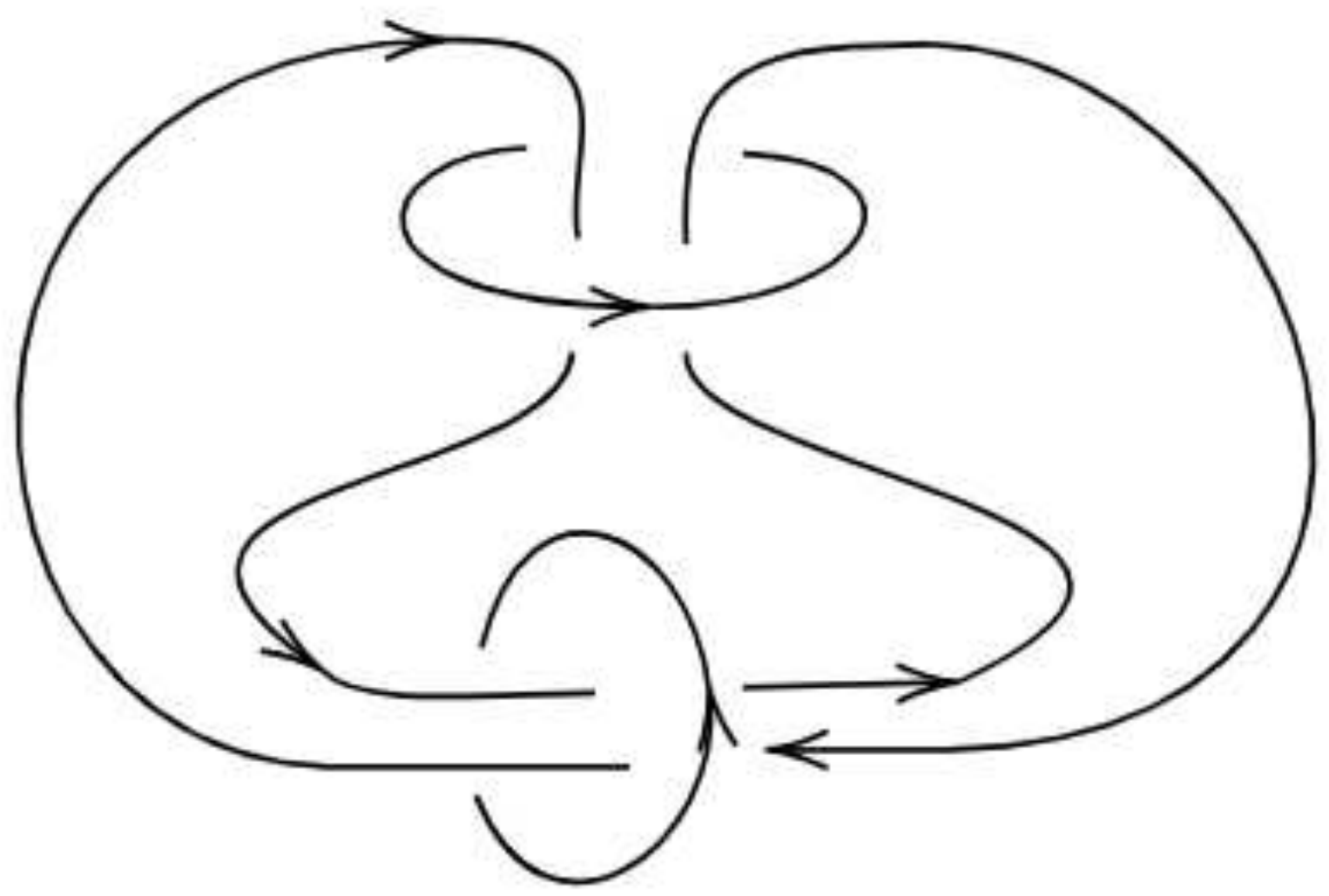}}
\end{align*}
\caption{Borromean ring $\mathbf{B}$} \label{fig:Borromean}
\end{figure}

The Borromean ring is a 3-component link shown in Figure \ref{fig:Borromean}, whose colored HOMFLYPT invariants is given by 

\begin{align}  \label{formulaB}
H_{r_1,r_2,r_3}(\mathbf{B};q,a)
=&\sum_{i=0}^{\min{\{r_1,r_2,r_3\}}}(-1)^i \frac{ \{ r_1+i-1;a
\}_{r_1-i} }{ \{ r_1-i\}! }\frac{ \{ r_2+i-1;a \}_{r_2-i} }{ \{
r_2-i\}! }\\\nonumber &\cdot\frac{ \{ r_3+i-1;a \}_{r_3-i} }{ \{
r_3-i\}! }\{2i-1;a\}_{2i}\{i-2;a\}_i\{i\}!. 
\end{align}
then
\begin{align}
\mathcal{H}_{r_1,r_2,r_3}(\mathbf{B}^{(\tau_1,\tau_2,\tau_3)};q,a)=(-1)^{\sum_{\alpha=1}^3r_\alpha\tau_\alpha}q^{\frac{\sum_{\alpha=1}^3r_\alpha(r_\alpha-1)\tau_\alpha}{2}}H_{r_1,r_2,r_3}(\mathbf{B}^{(\tau_1,\tau_2,\tau_3)};q,a).     
\end{align}

By straightforward computations, we obtain
\begin{align*}
p_{(1,1,1)}(\mathbf{B}^{(\tau_1,\tau_2,\tau_3)},1,a)&=-(-1)^{\sum_{\alpha=1}^3\tau_\alpha}a^{\frac{3}{2}}+3(-1)^{\sum_{\alpha=1}^3\tau_\alpha}a^{\frac{1}{2}}\\\nonumber
&-3(-1)^{\sum_{\alpha=1}^3\tau_\alpha}a^{-\frac{1}{2}}+(-1)^{\sum_{\alpha=1}^3\tau_\alpha}a^{-\frac{3}{2}},       
\end{align*}

\begin{align*}
&p_{2,1,1}(\mathbf{B}^{(\tau_1,\tau_2,\tau_3)},1,a)\\\nonumber
&=(-1)^{\tau_2+\tau_3}\left(-(1+\tau_1)a^2+(1-\tau_1)a^{-2}+(4\tau_1+2)a+(4\tau_1-2)a^{-1}-6\tau_1\right).    
\end{align*}

These formulas confirm the integrality of the BPS invariants for framed Borromean rings $\mathbf{B}^{(\tau_1,\tau_2,\tau_3)}$. 

Then, we list some Ooguri-Vafa invariants for Borromean rings in the following tables, which confirm the Conjecture \ref{conjecture-framedlink}.

\begin{table}[H]
\makebox[\textwidth]{
\centering
\renewcommand{\arraystretch}{1.25}
$\begin{array}{|c|cccc|}
\hline
i \backslash j & \frac{3}{2} & \frac{1}{2} & -\frac{1}{2} & -\frac{3}{2} \\ \hline
2 & 0 & -1 & 0 & 0 \\
1 & 1 & 1 & 0 & 0 \\
0 & -1 & 0 & 0 & 1 \\
-1 & 0 & 0 & -1 & -1 \\
-2 & 0 & 0 & 1 & 0 \\ \hline
\end{array}$
}
\newline \\ 
\caption{Ooguri-Vafa invariants $N_{(1,1,2),i,j}(\mathbf{B}^{(0,0,0)})$}
\label{}
\end{table}

\begin{table}
\makebox[\textwidth]{\centering
\renewcommand{\arraystretch}{1.25}
$\begin{array}{|c|cccc|}
\hline
i \backslash j & \frac{5}{2} & \frac{3}{2} & \frac{1}{2} & -\frac{1}{2} \\ \hline
2 & 0 & -1 & -1 & 0 \\
1 & 1 & 2 & 2 & 1 \\
0 & -1 & -2 & -2 & -1 \\
-1 & 0 & 1 & 1 & 0 \\ \hline
\end{array}$
}
\newline \\ 
\caption{Ooguri-Vafa invariants  $N_{(1,1,2),i,j}(\mathbf{B}^{(1,1,1)})$}
\label{}
\end{table}

\begin{table}
\makebox[\textwidth]{\centering
\renewcommand{\arraystretch}{1.25}
$\begin{array}{|c|ccccc|}
\hline
i \backslash j & 2 & 1 & 0 & -1 & -2 \\ \hline
\frac{5}{2} & 0 & -1 & 0 & 0 & 0 \\
\frac{3}{2} & 1 & 1 & -1 & 0 & 0 \\
\frac{1}{2} & -1 & 1 & 1 & 1 & 0 \\
-\frac{1}{2} & 0 & -1 & -1 & -1 & 1 \\
-\frac{3}{2} & 0 & 0 & 1 & -1 & -1 \\
-\frac{5}{2} & 0 & 0 & 0 & 1 & 0 \\ \hline
\end{array}$
}
\newline \\ 
\caption{Ooguri-Vafa invariants $N_{(1,2,2),i,j}(\mathbf{B}^{(0,0,0)})$}
\label{}
\end{table}

\begin{table}
\makebox[\textwidth]{\centering
\renewcommand{\arraystretch}{1.25}
$\begin{array}{|c|ccccc|}
\hline
i \backslash j & 4 & 3 & 2 & 1 & 0 \\ \hline
\frac{5}{2} & 0 & 1 & 2 & 1 & 0 \\
\frac{3}{2} & -1 & -3 & -4 & -3 & -1 \\
\frac{1}{2} & 1 & 3 & 4 & 3 & 1 \\
-\frac{1}{2} & 0 & -1 & -2 & -1 & 0 \\ \hline
\end{array}$
}
\newline \\ 
\caption{Ooguri-Vafa invariants $N_{(1,2,2),i,j}(\mathbf{B}^{(1,1,1)})$}
\label{}
\end{table}

\begin{table}
\makebox[\textwidth]{\centering
\renewcommand{\arraystretch}{1.25}
$\begin{array}{|c|cccccccccc|}
\hline
i \backslash j & \frac{9}{2} & \frac{7}{2} & \frac{5}{2} & \frac{3}{2} & \frac{1}{2} & -\frac{1}{2} & -\frac{3}{2} & -\frac{5}{2} & -\frac{7}{2} &
-\frac{9}{2} \\ \hline
3 & 0 & 1 & 0 & -2 & 0 & 1 & 0 & 0 & 0 & 0 \\
2 & -1 & -2 & 1 & 3 & -2 & 1 & 0 & 0 & 0 & 0 \\
1 & 2 & 2 & -1 & 1 & 1 & -2 & -1 & -1 & -1 & 0 \\
0 & -1 & -2 & -1 & -3 & 1 & -1 & 3 & 1 & 2 & 1 \\
-1 & 0 & 1 & 1 & 1 & 2 & -1 & -1 & 1 & -2 & -2 \\
-2 & 0 & 0 & 0 & 0 & -1 & 2 & -3 & -1 & 2 & 1 \\
-3 & 0 & 0 & 0 & 0 & -1 & 0 & 2 & 0 & -1 & 0 \\ \hline
\end{array}$
}
\newline \\ 
\caption{Ooguri-Vafa invariants  $N_{(2,2,2),i,j}(\mathbf{B}^{(0,0,0)})$}
\label{}
\end{table}

\begin{table}
\makebox[\textwidth]{\centering
\renewcommand{\arraystretch}{1.25}
$\begin{array}{|c|ccccccccccc|}
\hline
i \backslash j & \frac{15}{2} & \frac{13}{2} & \frac{11}{2} & \frac{9}{2} & \frac{7}{2} & \frac{5}{2} & \frac{3}{2} & \frac{1}{2} & -\frac{1}{2} &
-\frac{3}{2} & -\frac{5}{2} \\ \hline
3 & 0 & 1 & 0 & -2 & -3 & -3 & -1 & -1 & 0 & 0 & 0 \\
2 & -1 & -2 & 1 & 5 & 7 & 8 & 6 & 1 & 0 & -1 & 0 \\
1 & 2 & 2 & 0 & -3 & -7 & -9 & -7 & -2 & 0 & 2 & 1 \\
0 & -1 & -2 & -2 & -1 & 2 & 5 & 4 & 4 & 0 & 0 & -1 \\
-1 & 0 & 1 & 1 & 2 & 2 & -1 & -3 & -2 & -2 & -1 & 0 \\
-2 & 0 & 0 & 0 & -1 & -1 & 0 & 0 & 0 & 2 & 0 & 0 \\
-3 & 0 & 0 & 0 & 0 & 0 & 0 & 1 & 0 & 0 & 0 & 0 \\ \hline
\end{array}$
}
\newline \\ 
\caption{Ooguri-Vafa invariants $N_{(2,2,2),i,j}(\mathbf{B}^{(1,1,1)})$}
\label{}
\end{table}

\section{Discussions and further questions}
In this final section, we provide some related questions that deserve further study.

1. Extend the relationship between BPS invariants and the dual $A$-polynomial of framed knots to the cases of framed links. For a framed knot, there is an associated algebraic curve, named $a$-deformed $A$-polynomial, from which we obtain the dual $A$-polynomial. Then the BPS invariants of the framed knot can be read from this dual $A$-polynomial, see formula (\ref{formula-bpsfromcurve}) in Section \ref{section-dualAextremal}. As to a framed link, there is an associated higher-dimensional variety generalizing algebraic curve, hence, a natural question is how to extract the BPS invariants of a framed link from its associated variety? For the Whitehead link and Borromean rings, we have calculated a lot of BPS invariants for them. So they are the basic examples to be studied further in this direction. 

2. Calculate the Ooguri-Vafa invariants from the dual $A$-polynomials of the framed knots. It is known that the BPS invariant studied in this article is related to the open Gromov-Witten invariant in genus zero, see formula (\ref{formula-bps-GW}). The general Ooguri-Vafa invariants are related to the open Gromov-Witten invariants in higher genus that can be calculated using the topological recursions method developed in \cite{EO07,BKMP09,EO15,FLZ20}.  It is natural to consider applying the topological recursion method to the dual $A$-polynomial of the framed knot to calculate the corresponding Ooguri-Vafa invariants; see \cite{BEM12} for the application of the topological recursion method to torus knots.   

3. A challenge problem is how to prove Conjecture \ref{conjecture-framedlink}. Note that the statement of Conjecture \ref{conjecture-framedlink} is very general; it implies the integrality of BPS invariants, even for extremal BPS invariants given by formulas (\ref{formula-bn-p1}), (\ref{formula-bn-p2}), (\ref{formula-bn+p1}) and (\ref{formula-bn+p2}), the proof of integrality for them is nontrivial, see \cite{LZhu19}.  In \cite{CLPZ23}, we investigate the essential structures that appear in the HOMFLYPT skein theory to study the LMOV conjecture for framed links. From our experience, Conjecture \ref{conjecture-framedlink} could be easier to deal with than LMOV conjecture. For example, Conjecture \ref{conjecture-framedlink} has been checked/proved for some unframed knots/links; see \cite{ZR12,KRSS17,PSS18}.

\end{document}